\documentclass{emsprocart}

\usepackage{srcltx}
\usepackage{stmaryrd}
\usepackage[all]{xy}

\contact[christof@math.unam.mx]{1. Christof Geiss, 
Instituto de Matem\'aticas, Universidad Nacional Aut\'onoma de M\'exico, 
Ciudad Universitaria, 04510 M\'exico D. F. M\'exico}

\contact[bernard.leclerc@unicaen.fr]{2. Bernard Leclerc, 
LMNO UMR 6139, Universit\'e de Caen Basse-Normandie, CNRS, Campus 2, F-14032 Caen cedex, France.
Institut Universitaire de France}

\contact[schroer@math.uni-bonn.de]{3. Jan Schr\"oer, Mathematisches Institut, Universit\"at Bonn, Endenicher Allee 60, D-53115 Bonn, Germany}

\newtheorem{theorem}{Theorem}[section]

\newtheorem{prop}[theorem]{Proposition}
\newtheorem{conjecture}[theorem]{Conjecture}

\theoremstyle{definition}
\newtheorem{definition}[theorem]{Definition}
\newtheorem{remark}[theorem]{Remark}
\newtheorem{example}[theorem]{Example}

\newcommand{\Hom}{\operatorname{Hom}}

\newcommand{\End}{\operatorname{End}}
\newcommand{\Gr}{\operatorname{Gr}}

\def\N{{\mathbb N}}
\def\Z{{\mathbb Z}}
\def\C{{\mathbb C}}
\def\Q{{\mathbb Q}}
\def\P{{\mathbb P}}
\def\CC{\mathcal{C}}

\def\M{\mathcal{M}}
\def\A{\mathcal{A}}
\def\F{\mathcal{F}}

\def\eg{\emph{e.g.\,}}
\def\ie{\emph{i.e.\,}}
\def\resp{\emph{resp.\,}}

\def\ga{\gamma}

\def\vpi{\varpi}

\def\la{\lambda}
\def\Si{\Sigma}
\def\L{\Lambda}

\def\g{\mathfrak{g}}
\def\h{\mathfrak{h}}
\def\n{\mathfrak{n}}

\def\De{\Delta}

\def\<{\langle\,}
\def\>{\,\rangle}

\def\y{\mathbf{y}}

\def\SL{\mathrm{SL}}
\def\GL{\mathrm{GL}}
\def\Sp{\mathrm{Sp}}
\def\Gr{\mathrm{Gr}}
\def\a{\alpha}
\def\b{\beta}
\def\l{\lambda}
\def\ga{\gamma}
\def\rk{\mathrm{rank}}
\def\ad{\mathrm{ad}}
\def\slchap{\widehat{\mathfrak{sl}}}
\def\nil{\mathrm{nil}}
\def\i{\mathbf{i}}
\def\f{\mathfrak{f}}
\def\d{\mathbf{d}}

\def\vph{\varphi}
\def\Ext{\mathrm{Ext}}
\def\ZZ{\mathcal{Z}}
\def\SC{\mathcal{S}}

\def\stCC{\underline{\CC}}
\def\add{\mathrm{add}}
\def\G{\Gamma}

\def\j{\mathbf{j}}
\def\RR{\mathcal{R}}
\def\t{\mathbf{t}}
\def\sx{\mathbf{x}}
\def\T{\mathcal{T}} 
\def\Sub{\mathrm{Sub}}
\def\md{\mathrm{mod}}

\title[Cluster algebras in Lie theory]{Cluster algebras in algebraic Lie theory}

\author[Christof Geiss, Bernard Leclerc, Jan Schr\"oer]{Christof Geiss, Bernard Leclerc, Jan Schr\"oer}

\begin{document}

\maketitle

\begin{abstract}
We survey some recent constructions of cluster algebra structures on 
coordinate rings of unipotent subgroups and unipotent cells of Kac-Moody groups. 
We also review a quantized version of these results.
\end{abstract}



\maketitle

\bigskip
\setcounter{tocdepth}{1}
\tableofcontents

\section{Introduction}

Cluster algebras were invented by Fomin and Zelevinsky \cite{FZ1}
as an abstraction of certain combinatorial structures which they 
had previously discovered while studying total positivity in
semisimple algebraic groups.
A cluster algebra is a commutative ring with a distinguished set
of generators and a particular type of relations. 
Although there can be infinitely many generators and relations,
they are all obtained from a finite number of them
by means of an inductive procedure called mutation.
The precise definition of a cluster algebra will be recalled 
in \S\ref{sect2} below.

Several examples arising 
in Lie theory were already mentioned in \cite{FZ1},
like $\C[\SL_4/N]$, $\C[\Sp_4/N]$, $\C[\Gr_{2,n+3}]$.
Moreover, Fomin and Zelevinsky \cite[p.~498]{FZ1} conjectured that 
\begin{quote}
\emph{``the above examples can be extensively generalized:
for any simply-connected connected semisimple group $G$,
the coordinate rings $\C[G]$ and $\C[G/N]$, as well as coordinate
rings of many other interesting varieties related to $G$, have
a natural structure of a cluster algebra. This structure should
serve as an algebraic framework for the study of ``dual canonical
bases'' in these coordinate rings and their $q$-deformations.''} 
\end{quote}
In this survey we shall review the results of a series of papers 
\cite{GLS2, GLSUni1, GLSKM, GLSquantum},
in which we have implemented (part of) this program for 
certain varieties associated with elements $w$ of the Weyl group 
of a Kac-Moody group, namely, 
the unipotent subgroups $N(w)$ and the unipotent cells $N^w$.

\section{Cluster algebras}\label{sect2}

We start with a quick (and very incomplete) introduction to the theory
of cluster algebras. For more developed accounts we refer to
\cite{FominSurvey,FZSurvey,GSW,Keller}.

\subsection{Definitions and basic results}
Consider the field of rational functions $\F=\Q(x_1,\ldots,x_n)$.
A \emph{seed} in $\F$ is a pair $\Si=(\y,Q)$, where
$\y=(y_1,\ldots,y_n)$ is
a free generating set of $\F$, and 
$Q$ is a quiver (\ie an oriented graph) with vertices
labelled by $\{1,\ldots,n\}$.
We assume that $Q$ has neither loops nor 2-cycles.
For $k=1,\ldots,n$, one defines a new seed 
$\mu_k(\Si)=(\mu_k(\y),\mu_k(Q))$
as follows. First $\mu_k(y_i) = y_i$ for $i\not = k$,
and
\begin{equation}\label{FZmutation_form}
\mu_k(y_k)= \frac{\prod_{i\to k}y_i + \prod_{k\to j} y_j}{y_k}, 
\end{equation}
where the first (\resp second) product in the right hand side is over all arrows of 
$Q$ with target (\resp source) $k$.
Next $\mu_k(Q)$ is obtained from $Q$ by
\begin{itemize}
 \item[(a)] adding a new arrow $i\to j$ for every existing pair of arrows
$i\to k$ and $k\to j$;
 \item[(b)] reversing the orientation of every arrow with target or source
equal to $k$;
 \item[(c)] erasing every pair of opposite arrows possibly created by (a).
\end{itemize}
It is easy to check that $\mu_k(\Si)$ is a seed, and $\mu_k(\mu_k(\Si)) = \Si$.
The \emph{mutation class} $\CC(\Si)$ is the set of all seeds obtained 
from $\Si$ by a finite sequence of mutations $\mu_k$. One can think of
the elements of $\CC(\Si)$ as the vertices of an $n$-regular
tree in which every edge stands for a mutation.
If $\Si' = ((y'_1,\ldots,y'_n),Q')$ is a seed in $\CC(\Si)$, 
then the subset $\{y'_1,\ldots,y'_n\}$
is called a \emph{cluster}, and its elements are called \emph{cluster variables}.
 
Now, Fomin and Zelevinsky define the \emph{cluster algebra} $\A_\Si$ as the
subring of $\F$ generated by all cluster variables.
Some important elements of $\A_\Si$ are the \emph{cluster monomials},
\ie monomials in the cluster variables supported on a single cluster.

\begin{example}\label{ex2.1}
If $n=2$ and $\Si=((x_1,x_2),Q)$, where $Q$ is the quiver with $a$ 
arrows from $1$ to $2$,
then $\A_\Si$ is the subring of $\Q(x_1,x_2)$ generated 
by the rational functions $x_k$ defined recursively by 
\begin{equation}
x_{k+1}x_{k-1} = 1+ x_k^a,\qquad (k\in \Z).
\end{equation}
The clusters of $\A_\Si$ are the subsets $\{x_k,x_{k+1}\}$,
and the cluster monomials are the special elements of the form
\[
 x_k^l x_{k+1}^m,\qquad (k\in\Z,\ l,m\in \N).
\]
It turns out that when $a=1$, there are only five different clusters
and cluster variables, namely
\[
 x_{5k+1}=x_1,\ \
 x_{5k+2}=x_2,\ \
 x_{5k+3}=\frac{1+x_2}{x_1},\ \
 x_{5k+4}=\frac{1+x_1+x_2}{x_1x_2},\ \
 x_{5k}=\frac{1+x_1}{x_2}.
\]
For $a\ge 2$ though, the sequence $(x_k)$ is no longer periodic and $\A_\Si$
has infinitely many cluster variables. 
\end{example}
 
The next theorem summarizes the first deep results of this theory obtained by Fomin and Zelevinsky.
\begin{theorem}[\cite{FZ1},\cite{FZ2}]\label{thFZ}
\begin{itemize}
\item[\rm (i)] Every cluster variable of $\A_\Si$ is a Laurent polynomial 
with coefficients in $\Z$
in the cluster variables of any single fixed cluster.
\item[\rm (ii)] $\A_\Si$ has finitely many clusters if and only if the mutation
class $\CC(\Si)$ contains a seed whose quiver is an orientation of
a Dynkin diagram of type $A$, $D$, $E$.
\end{itemize}
\end{theorem}

\subsection{An example from Lie theory}\label{sect2.2}
We illustrate Theorem~\ref{thFZ}~(ii) with a prototypical example.
Let $G=SL_4$ and denote by $N$ the subgroup of upper unitriangular
matrices. 
In \cite[\S 2.6]{BFZ} explicit initial seeds for a 
cluster algebra structure in the coordinate ring of the big
cell of the base affine space $G/N$ were described. A simple modification
yields initial seeds for $\C[N]$ (see \cite{GLS3}). 
One of these seeds is
\[
((D_{1,2},\,D_{1,3},\,D_{12,23},\,D_{1,4},\,D_{12,34},\,D_{123,234}), Q), 
\]
where $Q$ is the triangular quiver:
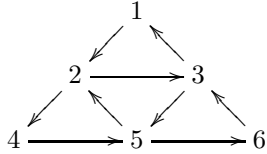
\[
\xymatrix@-1.0pc{
&& \ar[dl]{1} && \\
&\ar[dl]{2}\ar[rr]&&{3}\ar[ul]\ar[dl]\\
{4}\ar[rr]&&{5}\ar[ul]\ar[rr]&&{6}\ar[ul]
} 
\]
Here, by $D_{I,J}$ we mean the regular function on $N$ 
assigning to a matrix its minor with row-set $I$
and column-set $J$.
Moreover, the variables 
\[
x_4=D_{1,4},\ x_5=D_{12,34},\ x_6=D_{123,234}
\]
are \emph{frozen}, \ie they cannot be mutated, and therefore they belong
to every cluster. 
If one performs the mutation $\mu_1$, one gets the new cluster variable
\[
\mu_1(x_1) = \frac{x_2+x_3}{x_1} = \frac{D_{1,3}+D_{12,23}}{D_{1,2}} = D_{2,3},
\]
and thus the new seed
\[
((D_{2,3},\,D_{1,3},\,D_{12,23},\,D_{1,4},\,D_{12,34},\,D_{123,234}), \mu_1(Q)), 
\]
where $\mu_1(Q)$ is the mutated quiver:
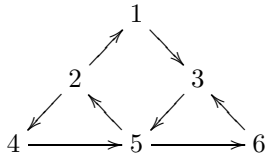
\[
\xymatrix@-1.0pc{
&& \ar[dr]{1} && \\
&\ar[dl]{2}\ar[ur]&&{3}\ar[dl]\\
{4}\ar[rr]&&{5}\ar[ul]\ar[rr]&&{6}\ar[ul]
} 
\]
The full subquiver of $\mu_1(Q)$ obtained by erasing vertices $4,5,6$
corresponding to the frozen variables is a Dynkin quiver of type $A_3$, hence,
by Theorem~\ref{thFZ}, this cluster algebra has finitely many clusters 
and cluster variables. It is an easy exercise to check that indeed it  
has $14$ clusters and $12$ cluster variables if we count the $3$ frozen ones.
Moreover, it follows from \cite{BZ0} that Lusztig's dual canonical basis of $\C[N]$ 
coincides with the set of all cluster monomials.


Finally, note that the open subset of $N$ given by the non-vanishing
of the 3 frozen variables $x_4, x_5, x_6$ is equal to 
\[
N^{w_0} := N \cap (B_-w_0B_-), 
\]
where $B_-$ denotes the subgroup of lower triangular matrices in $G$,
and $w_0$ is the longest element of the Weyl group of $G$.
This is an example of a \emph{unipotent cell}, that is, a stratum
of the decomposition of the unipotent group $N$ obtained by intersecting
it with the Bruhat decomposition of $G$ associated with the opposite
Borel subgroup $B_-$. The coordinate ring $\C[N^{w_0}]$
is obtained by localizing the polynomial ring $\C[N]$
at the element $x_4x_5x_6$, and we can see that it carries almost the same
cluster algebra structure as $\C[N]$, the only difference being that
the \emph{coefficient ring} generated by the frozen cluster variables
is now the Laurent polynomial ring in $x_4, x_5, x_6$.  

\section{Lie theory}

We now introduce a class of varieties generalizing the varieties 
$N$ and $N^{w_0}$ of Example~\ref{sect2.2}. These are subvarieties
of unipotent subgroups of symmetric Kac-Moody groups, and we 
first prepare the necessary notation and definitions. For more
details, we refer to \cite{Ku}.

\subsection{Kac-Moody Lie algebras}\label{sect3.1}
Let $C = [c_{ij}]$ be a symmetric  $n\times n$ generalized Cartan matrix.
Thus $2I - C$ is the adjacency matrix of an unoriented graph $\Gamma$ with $n$ vertices
and no loop.
Without loss of generality, we can assume that this graph is 
connected.
We will often denote by $I = [1,n]$ the indexing set of rows and
columns of $C$.
Let $\h$ be a $\C$-vector space of dimension $2n-\rk(C)$.
We choose linear independent subsets $\{h_i\mid i\in I\} \subset \h$
and $\{\a_i\mid i\in I\} \subset \h^*$ such that $\a_i(h_j) = c_{ij}$.
Let $\g$ be the Kac-Moody Lie algebra over $\C$ with generators
$e_i, f_i\ (i\in I)$, $h\in\h$, subject to the following relations:
\[
\begin{array}{l}
[h,h']=0,\ [h,e_i]={\a_i}(h)e_i,\ [h_,f_i]=-{\a_i}(h) f_i,\qquad (h,h'\in\h),\\[2mm]
[e_i,f_j]=\delta_{ij}h_i,  \hskip 5.7cm (i,j\in I),\\[3mm]
\ad(e_i)^{1-{c_{ij}}}(e_j) = 
\ad(f_i)^{1-{c_{ij}}}(f_j) = 0,\ \hskip 2.1cm(i\not = j).
\end{array}
\]
We denote by $\n_+$ (\resp $\n_-$) the subalgebra generated by $e_i\ (i\in I)$
(\resp $f_i\ (i\in I)$). For simplicity we shall often write $\n$ instead of $\n_+$.

Let $W$ be the subgroup of $\GL(\h^*)$ generated by the reflexions 
\[
s_i(\a) = \a - \a(h_i)\a_i,\qquad (i\in I,\ \a\in\h^*). 
\]
This is a Coxeter group with length function $w \mapsto \ell(w)$.
For $\a\in\h^*$ let 
\[
\g_\a = \{x\in \g \mid [h,x] = \a(h)x,\ h\in\h\}.
\]
We denote by $\De := \{\a\in\h^*\mid \g_\a \not = 0\}$ the root system of $\g$,
by $R:= \oplus_{i\in I} \Z\a_i$ the root lattice, by $R^+:= \oplus_{i\in I} \N\a_i$  
its positive cone, and by $\De^+ := \De\cap R^+$ the subset of positive roots.
We have $\De = \De^+\sqcup (-\De^+)$.
The Weyl group $W$ acts on $\De$, and we define the subset of real roots as the
$W$-orbit of $\{\a_i \mid i\in I\}$.
For $w\in W$, put
\[
 \De_w := \{\a\in\De^+ \mid w(\a) \in \De^-\}.
\]
This is a finite set of positive real roots, with cardinality $\ell(w)$.
Finally, set 
\[
\n({w}) := \bigoplus_{\a\in\De_{w}} \g_\a \subset \n_+,
\]
a nilpotent subalgebra of $\g$ of dimension $\ell(w)$. 

\begin{example}\label{example1}
Take 
\[
C = 
\begin{pmatrix}
2 & -2 \\
-2 & 2 
\end{pmatrix}.
\]
Then $\g=\slchap_2$ is an affine Lie algebra of Dynkin type $\widetilde{A}_1$.
We have 
\[
W = \langle s_1, s_2 \mid s_1^2 = s_2^2 = 1 \rangle.
\]
Let $w = s_2s_1s_2s_1$. Then
\[
\De_w = \{{\a_1},\, 2{\a_1} + {\a_2} ,\, 3{\a_1} + 2{\a_2} ,\, 4{\a_1} + 3{\a_2}\},
\]
and 
\[
\n(w) = \mathrm{Span}\langle e_1,\ [e_1,[e_2,e_1]], 
[e_1,[e_2,[e_1,[e_2,e_1]]]],\ [e_1,[e_2,[e_1,[e_2,[e_1,[e_2,e_1]]]]]]\rangle. 
\]
\end{example}

\subsection{Kac-Moody groups}

The enveloping algebra $U(\n)$ of the Lie algebra $\n$ is a cocommutative Hopf algebra, with
an $R^+$-grading given by $\deg(e_i) = \a_i$.
Let 
\[
U(\n)_{\mathrm{gr}}^* := \bigoplus_{d\in R^+} U(\n)_d^*
\]
be its graded dual.
This is a commutative Hopf algebra.
Define 
\[
N := \mathrm{max\,Spec}(U(\n)_{\mathrm{gr}}^*) = \Hom_{\mathrm{alg}}(U(\n)_{\mathrm{gr}}^*,\C).
\]
The comultiplication of $U(\n)_{\mathrm{gr}}^*$ gives $N$ a group structure.
As a group, $N$ can be identified with the pro-unipotent pro-group with Lie algebra 
\[
 \widehat{\n} = \prod_{{\a}\in\De^+} \g_{\a}.
\]
By construction, we can identify $U(\n)_{\mathrm{gr}}^*$ with the coordinate ring $\C[N]$
of $N$.

For $w\in W$, let $N(w)$ be the subgroup of $N$
with Lie algebra $\n(w)$.
Define also $N'(w)$ to be the subgroup of $N$
with Lie algebra 
\[
\n'(w) := \prod_{{\a}\not\in\De_{w}} {\g}_{\a} 
\subset \widehat{{\n}}.
\]
Multiplication yields a bijection 
$N({w})\times N'({w}) \stackrel{\sim}{\to} N$.
\begin{prop}[{\cite[Proposition 8.2]{GLSKM}}]\label{Prop3.2}
The coordinate ring $\C[{N}({w})]$ is isomorphic to the 
invariant subring 
\[
\C[{N}]^{{N}'({w})} = 
\left\{f\in\C[{N}] \mid f(nn') = f(n),\ n\in {N},\ n'\in N'(w)\right\}.
\]
\end{prop}
Let $G$ be the group attached to $\g$ by Kac and Peterson \cite{KP}.
This is an affine ind-variety.
It has a refined Tits system 
\[
({G}, \mathrm{Norm}_{{G}}({H}),{N_+},{N_-},{H}),
\]
where $\mathrm{Lie}({H})={\h}$, $\mathrm{Lie}({N_+})={\n_+}$, and $\mathrm{Lie}({N_-})={\n_-}$
(see \cite{Ku}).
Note that in general $N \not \subset G$. Both $N$ and $G$ can be regarded as subgroups of 
a bigger group ${G}^{\mathrm{max}}$ constructed by Tits. Then ${N_+} = {N} \cap {G}$.

\begin{example}\label{Example2}
(a) If $\g$ is finite-dimensional of Dynkin type $X_n$, then $G = G_{X_n}(\C)$ is a 
connected simply-connected algebraic group of type $X_n$ over $\C$, and $N_+$
is a maximal unipotent subgroup $N_+ = N_{X_n}(\C)$. 
Thus, if $X_n = A_n$, ${G} = \SL(n+1,\C)$, and $N_+$ is the subgroup of unipotent
upper triangular matrices.

(b) If $\g$ is an affine Lie algebra of affine Dynkin type $\widetilde{X}_n$, 
then $G$ is a central extension of $G_{X_n}(\C[z,z^{-1}])$ by $\C^*$. 
Moreover,
\[
{N_+} \simeq \{ g \in {G}_{X_n}(\C[z]) \mid g|_{z=0} \in {N}_{X_n}(\C) \}. 
\]
Thus, continuing Example~\ref{example1}, if $\g$ is of type $\widetilde{A}_1$,
\[
N_+ = \left\{
\begin{pmatrix}
a & b \\
c & d 
\end{pmatrix}
\in \SL(2,\C[z])\mid
 a(0)=d(0)=1,\ c(0)=0 
\right\}.  
\]

(c) If $\g$ is of indefinite type, no ``concrete'' realization of ${G}$ is known. 
\end{example}

\subsection{Generalized minors}\label{sect3.3}
 
We have $\mathrm{Norm}_{{G}}({H})/{H} \cong {W}$.
For $i\in{I}$, put 
\[
\overline{s}_i := \exp(f_i)\exp(-e_i)\exp(f_i),
\quad
\overline{\overline{s}}_i := \exp(-f_i)\exp(e_i)\exp(-f_i).
\]
To ${w} = {s_{i_1}}\cdots {s_{i_r}} \in W$ with $\ell({w})=r$, 
we attach two representatives in $\mathrm{Norm}_{{G}}({H})$:
\[
\overline{{w}}=\overline{s}_{i_1}\cdots \overline{s}_{i_r},
\qquad
\overline{\overline{w}}=\overline{\overline{s}}_{i_1}\cdots \overline{\overline{s}}_{i_r}.
\]
Let ${G_0} = {N_-} {H} {N_+}$ be the Zariski open subset of $G$ consisting of
elements $g$ having a Birkhoff decomposition. For $g\in G_0$, this unique 
factorization is written 
\[
g=[g]_-[g]_0[g]_+,\qquad
([g]_-\in{N_-},\ [g]_0\in{H},\  [g]_+\in{N_+}). 
\]
Let $\{\varpi_i \mid i\in I\} \subset \h^*$ be a fixed choice
of fundamental weights, that is, 
\[
\varpi_i(h_j) = \delta_{ij}, \qquad (i,j\in I).
\]
Let $x \mapsto x^{\vpi_i}$ denote the corresponding characters of ${H}$. 
There is a unique regular function $\De_{{\vpi_i},{\vpi_i}}$
on $G$ such that
\[
\De_{{\vpi_i},{\vpi_i}}(g) = [g]_0^{\vpi_i}, \qquad (g\in{G_0}).
\]
Moreover, 
 $G_0 = \{ g\in G \mid \De_{{\vpi_i},{\vpi_i}}(g) \not = 0,\ i\in{I}\}$.

\begin{definition}
For $u,v\in W$ and $i\in I$, the generalized minor
$\De_{u(\vpi_i),v(\vpi_i)}$ is the regular function on $G$
given by
\[
\De_{u(\vpi_i),v(\vpi_i)}(g) = \De_{{\vpi_i},{\vpi_i}}\left(\overline{\overline{u^{-1}}} g\overline{v}\right),
\qquad (g\in G).
\]
\end{definition}

\subsection{Unipotent cells}
Let ${B_-} = {N_- H}$.
The group ${G}$ has a Bruhat decomposition: 
\[
{G} = \bigsqcup_{{w}\in{W}} {B_-}\overline{{w}}{B_-}.
\]
For ${w}\in{W}$, define the \emph{unipotent cell}
${N}^{w}:= {N_+}\,\cap\,({B_-}\overline{w}{B_-})$.
Let 
\[
O_{w} := \{g\in {N}({w})\mid \De_{{\vpi_i},{w}^{-1}({\vpi_i})}(g) \not = 0,\ i\in{I}\},
\]
an open subset of the affine space $N(w)$.

\begin{prop}[{\cite[Proposition 8.5]{GLSKM}}]\label{Prop3.5}
We have an isomorphism $O_{{w}} \stackrel{\sim}{\to} {N}^{{w}}$.
It follows that $\C[{N}^{{w}}]$ is the localization of
$\C[{N}({w})] \simeq \C[{N}]^{N'({w})}$ at
$\prod_{i\in {I}} \De_{{\vpi_i},{w}^{-1}({\vpi_i})}$.
\end{prop}

\subsection{Formulas for factorization parameters}\label{sect3.5}
Define the one-parameter subgroups of $N_+$:
\[
x_{{i}}({t}) := \exp({t}e_i),\qquad ({t}\in\C,\, {i}\in{I}).
\]
Consider as above a reduced decomposition ${w} = {s_{i_r}}\cdots {s_{i_1}} \in W$.
(Note that for reasons which will become clear in \S\ref{sect4.3} below,
we prefer from now on to number the factors of this decomposition
\emph{from right to left}.)
Then the image of the map $\mathbf{x}_{i_r,\ldots,i_1}\colon(\C^*)^r \to {N}$ given by
\[
 \mathbf{x}_{i_r,\ldots,i_1}({t_r},\ldots,{t_1}) := x_{{i_r}}({t_r})\cdots x_{{i_1}}({t_1})
\]
is a dense subset of ${N}^{{w}}$. 
More precisely, $\mathbf{x}_{i_r,\ldots,i_1}$ gives a birational isomorphism
from $(\C^*)^r$ to ${N}^{{w}}$.
 
Thus, if $f\in\C[{N}^{{w}}]$ then 
$f(\mathbf{x}_{i_r,\ldots,i_1}(\mathbf{t}))$
is a rational function of $\mathbf{t}:=(t_r,\ldots,t_1)$ which completely determines $f$.
Conversely, following Berenstein, Fomin and Zelevinsky \cite{BFZ0}, one can give explicit formulas
in terms of generalized minors for expressing each $t_i$ as a rational function on $N^w$ \cite{GLSCA}.
These \emph{Chamber Ansatz formulas} will be explained in \S\ref{sectCA} below.

\begin{example}\label{Example3.6}
We continue Example~\ref{example1}. Thus $\g = \slchap_2$ and $w=s_2s_1s_2s_1$. 
By Example \ref{Example2} (b), we have
\[
N_+ = \left\{\begin{pmatrix}
a & b \\
c & d 
\end{pmatrix}
\mid
a, b, c, d \in \C[z],\ ad-bc=1,\ 
 a(0)=d(0)=1,\ c(0)=0 
\right\}.  
\]
Writing, for $x=
\begin{pmatrix}
a & b \\
c & d 
\end{pmatrix}\in N_+$,
\[
a=1+\sum_{k\ge 1}a_kz^k,\quad
b=\sum_{k\ge 0}b_kz^k,\quad
c=\sum_{k\ge 1}c_kz^k,\quad
d=1+\sum_{k\ge 1}d_kz^k,
\]
we get regular functions $a_k(x), b_k(x), c_k(x), d_k(x)$ on $N_+$,
which by restriction give regular functions on $N^w$.
We have
\[
x_1(t) =  \begin{pmatrix}
1 & t \\
0 & 1 
\end{pmatrix},\qquad
x_2(t) =  \begin{pmatrix}
1 & 0 \\
tz & 1 
\end{pmatrix},
\]
and, putting $\mathbf{x}_{2,1,2,1}(\mathbf{t}) := x_2(t_4)x_1(t_3)x_2(t_2)x_1(t_1)$, 
we calculate
\[
\mathbf{x}_{2,1,2,1}(\mathbf{t}) = 
\begin{pmatrix}
1 + t_2t_3z& t_1+t_3+t_1t_2t_3z \\
(t_2+t_4)z+t_2t_3t_4z^2 & 1 +(t_1t_2+t_1t_4+t_3t_4)z+t_1t_2t_3t_4z^2 
\end{pmatrix}. 
\]
Hence
\[
\begin{array}{ll}
a_1(\mathbf{x}_{2,1,2,1}(\mathbf{t})) = t_2t_3,\\[2mm]
b_0(\mathbf{x}_{2,1,2,1}(\mathbf{t})) = t_1+t_3, 
\quad
&b_1(\mathbf{x}_{2,1,2,1}(\mathbf{t})) = t_1t_2t_3, \\[2mm]
c_1(\mathbf{x}_{2,1,2,1}(\mathbf{t})) = t_2+t_4,
\quad
&c_2(\mathbf{x}_{2,1,2,1}(\mathbf{t})) = t_2t_3t_4, \\[2mm]
d_1(\mathbf{x}_{2,1,2,1}(\mathbf{t})) = t_1t_2+t_1t_4+t_3t_4, 
\quad
&d_2(\mathbf{x}_{2,1,2,1}(\mathbf{t})) = t_1t_2t_3t_4.
\end{array}
\]
It follows that, 
\begin{equation}\label{eqnotCA}
t_4 = \frac{c_2}{a_1},\qquad
t_3 = \frac{\left|
\begin{matrix}
b_0& b_1\\
1 & a_1
\end{matrix}
\right|}{a_1},\qquad
t_2 = \frac{\left|
\begin{matrix}
c_1& c_2\\
1 & a_1
\end{matrix}
\right|}{a_1},\qquad 
t_1= \frac{b_1}{a_1}.
\end{equation}
In fact, these formulas, although similar in spirit, are \emph{not} the Chamber
Ansatz formulas. 
Indeed the regular functions $a_1$, $b_1$, $c_2$, $b_0a_1-b_1$, and $c_1a_1-c_2$,
are \emph{not} cluster variables for the cluster algebra structure on $\C[N^w]$. 
Compare Example~\ref{example6.2} below.
\end{example}

\begin{example}\label{Example3.7}
For the sake of comparison with Example~\ref{Example3.6}, let us describe 
the groups $N$, $N(w)$, and $N'(w)$, for the same $w=s_2s_1s_2s_1$ and $\g = \slchap_2$.
First, we have 
\[
N = \left\{\begin{pmatrix}
a & b \\
c & d 
\end{pmatrix}
\mid
a, b, c, d \in \C[[z]],\ ad-bc=1,\ 
 a(0)=d(0)=1,\ c(0)=0 
\right\},  
\]
where $\C[[z]]$ is the ring of formal power series in $z$.
Next,
\[
\begin{array}{lcl}
N(w) &=&  \left\{\begin{pmatrix}
1 & b_0+b_1z+b_2z^2+b_3z^3 \\
0 & 1 
\end{pmatrix}
\mid
b_0, b_1, b_2, b_3 \in \C
\right\},\\  [5mm]
N'(w) &=&  \left\{\begin{pmatrix}
a & b \\
c & d 
\end{pmatrix} \in N
\mid
b \in z^4\C[[z]]
\right\}.
\end{array}
\]
Finally, for $g\in N(w)$, we have 
\[
\De_{\varpi_1,w^{-1}(\varpi_1)}(g) = b_0b_2-b_1^2,\qquad
\De_{\varpi_2,w^{-1}(\varpi_2)}(g) = b_1b_3-b_2^2, 
\]
so 
\[
N^w \cong \left\{\begin{pmatrix}
1 & b_0+b_1z+b_2z^2+b_3z^3 \\
0 & 1 
\end{pmatrix}
\mid b_i\in\C,\ 
(b_0b_2-b_1^2)(b_1b_3-b_2^2) \not = 0
\right\}. 
\]

\end{example}

\section{Categories of modules over preprojective algebras}

Seminal works of Ringel and Lusztig have shown that
the interaction between Kac-Moody algebras and the 
representation theory of quivers is essential for understanding 
the quantum enveloping algebra $U_q(\n)$ and its canonical basis.
Following Lusztig, one can also construct the (classical) enveloping algebra
$U(\n)$ in terms of the path algebra of a quiver with relations
called the preprojective algebra. This will be
our basic tool for exploring cluster algebra structures
on $\C[N(w)]$ and $\C[N^w]$. In this section we shall introduce 
the preprojective algebra $\L$ and its nilpotent representations,
and explain how it yields an interesting basis of $\C[N]$
dual to Lusztig's semicanonical basis of $U(\n)$. 
We will then describe some categories $\CC_w$ of $\L$-modules
introduced and studied in general by Buan, Iyama, Reiten and Scott \cite{BIRS}
(and independently in \cite{GLSUni1}, for adaptable Weyl group elements~$w$).
This will provide a categorical model for the cluster algebras of
\S\ref{sect5}.

\subsection{The preprojective algebra}
Let $Q$ be a quiver obtained by orienting the edges of the graph $\Gamma$
of \S\ref{sect3.1}. We require $Q$ to be acyclic, that is, $Q$ has no oriented
cycle. 
Let $\overline{Q}$ denote the double quiver obtained 
from $Q$ by adjoining to every arrow $a\colon i\to j$
an opposite arrow $a^*\colon j \to i$.
Consider the element
\[
 \rho = \sum (aa^* - a^*a)
\]
of the path algebra $\C\overline{Q}$ of $\overline{Q}$,
where the sum is over all arrows $a$ of $Q$.
Following \cite{GP,Ri}, we define the \emph{preprojective 
algebra} $\L$ as the quotient of $\C\overline{Q}$ by 
the two-sided ideal generated by $\rho$.
It is well-known that $\L$ is independent of the choice 
of orientation $Q$ of $\Gamma$.
Moreover, $\L$
is finite-dimensional if and only if 
the Kac-Moody algebra $\g$ is finite-dimensional, that is,
if and only if $\Gamma$ is a Dynkin diagram of type $A, D, E$.

We say that a finite-dimensional $\L$-module
is \emph{nilpotent} if all its composition factors 
are one-dimensional. Let $\nil(\L)$ denote the category
of nilpotent $\L$-modules. This is an abelian category
with infinitely many isomorphism classes of indecomposable
objects, except if $\g$ has type $A_n$ with $n\le 4$.
It is remarkable that these few exceptional cases 
coincide precisely with the cases when the cluster algebra
$\C[N]$ has finitely many cluster variables. Moreover, 
it is a nice exercise to verify that the number of indecomposable 
$\L$-modules is then equal to the number of cluster
variables.   
This suggests a close relationship in general between
$\L$ and $\C[N]$. 
To describe it we start with Lusztig's Lagrangian
construction of the enveloping algebra $U(\n)$ \cite{Lu0, Lu1}.
This is a realization of $U(\n)$ as an algebra of $\C$-valued
constructible functions over the varieties of nilpotent representations of $\L$.

Denote by $S_i \ (1\le i\le n)$ the one-dimensional $\L$-module
supported on the vertex $i$ of $\overline{Q}$. 
Given a sequence $\i=(i_1,\ldots,i_d)$ and a nilpotent $\L$-module $X$
of dimension~$d$,
we introduce the variety $\F_{X,\i}$ of
flags of submodules
\[
\f = (0=F_0\subset F_1 \subset \cdots \subset F_d = X) 
\]
such that $F_k/F_{k-1} \cong S_{i_k}$ for $k=1,\ldots,d$.
This is a projective variety.
Denote by $\L_\d$ the affine variety of nilpotent $\L$-modules $X$ with a given
dimension vector $\d=(d_i)$, where $\sum_i d_i =d$.
Consider the constructible function $\chi_\i$ on $\L_\d$  
given by 
\[
\chi_\i(X) = \chi(\F_{X,\i})
\]
where $\chi$ denotes the Euler-Poincar\'e characteristic.
Let $\M_{\d}$ be the $\C$-vector space spanned by
the functions $\chi_\i$ for all possible sequences $\i$
of length $d$, and let 
\[
\M = \bigoplus_{\d\in\N^n} \M_\d. 
\]
Lusztig has endowed $\M$ with an associative multiplication
which formally resembles a convolution product, and he has
shown that, if we denote by $e_i$ the Chevalley generators of $\n$,
there is an algebra isomorphism $U(\n) \stackrel{\sim}{\rightarrow} \M$
mapping the product $e_{i_1}\cdots e_{i_d}$ to $\chi_\i$ for every $\i=(i_1,\ldots,i_d)$.

Now, following \cite{GLS1,GLS2}, we dualize the picture. 
Every $X\in\nil(\L)$ determines a linear form $\delta_X$ on $\M$ given by
\[
 \delta_X(f) = f(X),\qquad (f\in\M).
\]
Through the isomorphisms $\M_{\mathrm{gr}}^* \simeq U(\n)_{\mathrm{gr}}^* \simeq \C[N]$, the 
form $\delta_X$ corresponds to an element $\vph_X$ of $\C[N]$,
and we have thus attached to every object $X$ in $\nil(\L)$
a polynomial function $\vph_X$ on $N$.
Unwrapping this definition, we have an explicit formula for evaluating
$\vph_X$ on an arbitrary product 
$\mathbf{x}_{\mathbf{j}}(\mathbf{t}) := x_{{j_1}}({t_1})\cdots x_{{j_r}}({t_r})$, namely
\begin{equation}
\vph_X(\mathbf{x}_{\mathbf{j}}(\mathbf{t})) 
= 
\sum_{\mathbf{a}\in\N^r} \chi_{{\mathbf{j}}^{\mathbf{a}}}(X) \frac{\mathbf{t}^{\mathbf{a}}}{\mathbf{a}!},
\end{equation}
where 
$\displaystyle\frac{\mathbf{t}^{\mathbf{a}}}{\mathbf{a}!} := \frac{t_1^{a_1}\cdots t_r^{a_r}}{a_1!\cdots a_r!}$,
and
${\mathbf{j}}^{\mathbf{a}} := (j_1,\ldots,j_1,j_2, \ldots,j_2, \ldots, j_r\ldots,j_r)$ with 
each component $j_k$ repeated $a_k$ times.

\begin{example}
If $\g$ is of type $A_3$, and if we denote by $P_i$
the projective cover of~$S_i$, one has 
\[
 \vph_{P_1} = D_{123,234},\quad
 \vph_{P_2} = D_{12,34},\quad
 \vph_{P_3} = D_{1,4}.
\]
More generally, the functions $\vph_X$ corresponding to the
12 indecomposable $\L$-modules are 
the 12 cluster variables of $\C[N]$ (see \S\ref{sect2.2}).
\end{example}

Via the correspondence $X \mapsto \vph_X$, the ring $\C[N]$
can be regarded as a kind of (dual) Hall algebra of the category 
$\nil(\L)$. Indeed, the multiplication of $\C[N]$ encodes
extensions in $\nil(\L)$, as shown by the following crucial
result. Before stating it, we recall that $\nil(\L)$ possesses
a remarkable symmetry with respect to extensions, namely,
$\Ext^1_\L(X,Y)$ is isomorphic to the dual of 
$\Ext^1_\L(Y,X)$ functorially in $X$ and~$Y$ (see \cite{GLS4}).
In particular $\dim\Ext^1_\L(X,Y) = \dim \Ext^1_\L(Y,X)$ for
every~$X, Y$.

\begin{theorem}[\cite{GLS1,GLS4}]
\label{Th_mult}
Let $X, Y \in \nil(\L)$. 
\begin{itemize}
 \item[\rm(i)] We have 
$\vph_X \vph_Y = \vph_{X\oplus Y}$.
 \item[\rm(ii)] Assume that $\dim \Ext^1_\L(X,Y)=1$, and let
\[
 0\to X \to L \to Y \to 0 
 \quad\mbox{and}\quad 
 0\to Y \to M \to X \to 0
\]
be non-split short exact sequences. Then
$\vph_X \vph_Y = \vph_L + \vph_M$. 
\end{itemize}
\end{theorem}

\begin{example}
We illustrate Theorem~\ref{Th_mult}~(ii) in type $A_2$. Take 
$X=S_1$ and $Y=S_2$. Then we have the non-split short
exact sequences
\[
 0\to S_1 \to P_2 \to S_2 \to 0 
 \quad\mbox{and}\quad 
 0\to S_2 \to P_1 \to S_1 \to 0,
\]
which imply the relation 
\[
\vph_{S_1} \vph_{S_2} = \vph_{P_2} + \vph_{P_1},
\]
that is, the elementary determinantal relation 
\[
D_{1,2}D_{2,3} = D_{1,3} + D_{12,23}
\]
on the unitriangular subgroup of $\SL(3,\C)$.
\end{example}

\subsection{The dual semicanonical basis}

We can now introduce the dual semicanonical basis of the vector space $\C[N]$.
Let $\d = (d_i)$ be a dimension vector. The variety $\overline{E}_\d$
of representations of $\C\overline{Q}$ with dimension vector $\d$ is a
vector space 
with a natural symplectic structure.
Lusztig \cite{Lu0} has shown that $\L_\d$ is a Lagrangian
subvariety of $\overline{E}_\d$,
whose number of irreducible components is equal to 
the dimension of the degree $\d$ homogeneous component of $U(\n)$.
Let $Z$ be an irreducible component of $\L_\d$. Since 
$\vph \colon X \mapsto \vph_X$
is a constructible map on $\L_\d$, it is constant on a Zariski dense open subset of $Z$.
Let $\vph_Z$ denote this generic value of $\vph$ on~$Z$.
Then, if we denote by $\ZZ$ the collection of all
irreducible components of all varieties~$\L_\d$, one can easily check
that
\[
 \SC^* := \{\vph_Z \mid Z \in \ZZ\} \subset \C[N] 
\]
is dual to the basis $\SC = \{f_Z \mid Z \in \ZZ\}$ of $\M \cong U(\n)$
constructed by Lusztig in \cite{Lu1}, and called by him the semicanonical basis.

For example, 
%
suppose that $X\in\nil(\L)$ is \emph{rigid}, 
\ie that $\Ext^1_\L(X,X) = 0$.
Then $X$ is a generic point of the unique irreducible 
component $Z$ on which it sits, that is, $\vph_X = \vph_Z$
belongs to $\SC^*$.

\subsection{Categories attached to Weyl group elements}\label{sect4.3}

Let $w\in W$ be of length~$r$ and choose a reduced decomposition
$w=s_{i_r}\cdots s_{i_1}$ (note that again, we number the factors from
right to left).
It is well-known that $\De_w$ consists of the roots
\[
\b_k:=s_{i_1}\cdots s_{i_{k-1}}(\a_{i_k}),\qquad  (1\le k \le r).
\]
The following sequence of $R^+$ will also play an important role:
\[
\ga_k:=\varpi_{i_k} - s_{i_1}\cdots s_{i_{k}}(\varpi_{i_k}),\qquad  (1\le k \le r). 
\]
From now on, we shall freely identify dimension vectors $\d = (d_i)$ with
elements of $R^+$ via
\[
\d \equiv \sum_{i\in I} d_i \a_i. 
\]
For $k = 1,\ldots, r$, one can show that there is a unique $V_k \in \nil(\L)$ (up to isomorphism)
whose socle is $S_{i_k}$ and whose dimension vector is $\ga_k$.
Let us write $\i = (i_r,\ldots,i_1)$, and define 
\[
 V_\i := \bigoplus_{k=1}^r V_k \in \nil(\L).
\]
Up to duality, this $\L$-module is the same as the one introduced
by Buan, Iyama, Reiten, and Scott \cite{BIRS}.
Following \cite{BIRS}, we consider
the full subcategory of $\nil(\L)$ whose objects are factor modules of 
direct sums of finitely many copies of~$V_\i$.
It turns out that this category depends only on $w$ (and not on the
choice of a reduced expression $\i$), so we may denote it by $\CC_w$.

For $j\in I$, let $k_j := \max\{1\le k\le r\mid i_k = j\}$.
The submodule
$\bigoplus_{j\in I} V_{k_j}$ of $V_\i$ also depends only on $w$,
and we denote it by $I_w$.

\begin{example}\label{Example4.4}
We consider again $\g = \slchap_2$ and $w=s_2s_1s_2s_1$, so $\i=(2,1,2,1)$.
We take for $Q$ the Kronecker quiver:
$$
\xymatrix{1 \ar@<0.5ex>[r]^{a}\ar@<-0.5ex>[r]_b & 2}
$$
The following pictures describe the indecomposable direct summands
of $V_\i$. The numbers 1 and 2 in the pictures are basis vectors of the modules $V_k$.
The solid edges show how the arrows $a$ and $b$ of $\overline{Q}$ act on these
vectors, and the dotted edges illustrate the actions of $a^*$ and $b^*$.
The arrows $a$ and $b^*$ are pointing south east, and the arrow $a^*$ and $b$
are pointing south west.
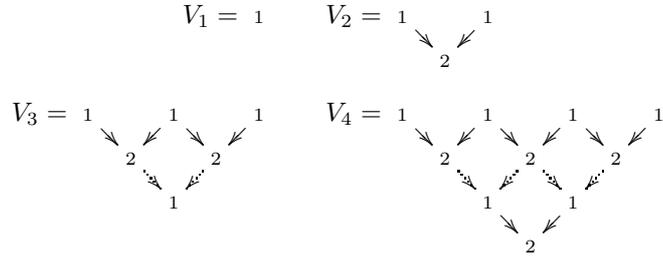
\[
\def\objectstyle{\scriptstyle}\def\labelstyle{\scriptstyle}
V_1=
\xymatrix@-1.5pc{
1
}
\qquad
\def\objectstyle{\scriptstyle}\def\labelstyle{\scriptstyle}
V_2=
\xymatrix@-1.5pc{
1\ar@{->}[rd]& &\ar@{->}[ld]1 \\
 & 2
}
\] 
\[
\def\objectstyle{\scriptstyle}\def\labelstyle{\scriptstyle}
V_3=
\xymatrix@-1.5pc{
 1\ar@{->}[rd]& &\ar@{->}[ld]1\ar@{->}[rd]& &\ar@{->}[ld]1\\
  &2\ar@{.>}[rd]& &\ar@{.>}[ld]2 \\
 & & 1
}
\qquad
\def\objectstyle{\scriptstyle}\def\labelstyle{\scriptstyle}
V_4=
\xymatrix@-1.5pc{
1\ar@{->}[rd]&&\ar@{->}[ld]1\ar@{->}[rd]&&\ar@{->}[ld]1\ar@{->}[rd]&&\ar@{->}[ld]1\\
 &2\ar@{.>}[rd]& &\ar@{.>}[ld]2\ar@{.>}[rd]& &\ar@{.>}[ld]2\\
 & &1\ar@{->}[rd]& &\ar@{->}[ld]1 \\
& & & 2
}
\] 
Here, $I_w = V_3\oplus V_4$. The modules
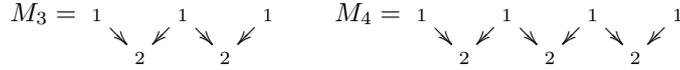
\[
\def\objectstyle{\scriptstyle}\def\labelstyle{\scriptstyle}
M_3=
\xymatrix@-1.5pc{
 1\ar@{->}[rd]& &\ar@{->}[ld]1\ar@{->}[rd]& &\ar@{->}[ld]1\\
  &2& &2
}
\qquad
\def\objectstyle{\scriptstyle}\def\labelstyle{\scriptstyle}
M_4=
\xymatrix@-1.5pc{
1\ar@{->}[rd]&&\ar@{->}[ld]1\ar@{->}[rd]&&\ar@{->}[ld]1\ar@{->}[rd]&&\ar@{->}[ld]1\\
 &2& &2& &2\\
}
\] 
are factor modules of $V_3$ and $V_4$ respectively, so they are objects of $\CC_w$.
On the other hand the $\L$-module
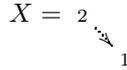
\[
\def\objectstyle{\scriptstyle}\def\labelstyle{\scriptstyle}
X=
\xymatrix@-1.5pc{
2\ar@{.>}[rd]&  \\
 & 1
}
\]
cannot be obtained as a factor module of (a direct sum of copies) of $V_\i$,
hence it does not belong to $\CC_w$.

\end{example}

\subsection{Frobenius categories}
Let us recall some definitions from homological algebra.
Let $\CC$ be a subcategory of the category of modules
over an algebra $A$,  
which is closed under extensions.
Clearly, we have
$$
\Ext_{\CC}^1(X,Y) = \Ext_A^1(X,Y)
$$ 
for all modules $X$ and $Y$ in $\CC$.
An $A$-module $C$ in $\CC$ is called $\CC$-{\it projective}  
(resp. $\CC$-{\it injective})
if $\Ext_A^1(C,X) = 0$ (resp. $\Ext_A^1(X,C) = 0$) for all $X \in \CC$.
If $C$ is $\CC$-projective and $\CC$-injective, then $C$ is also called
$\CC$-{\it projective-injective}.
We say that $\CC$ has {\it enough projectives} 
(resp. {\it enough injectives})
if for each $X \in \CC$ there exists a short exact sequence
$0 \to Y \to C \to X \to 0$
(resp. $0 \to X \to C \to Y \to 0$)
where $C$ is $\CC$-projective (resp. $\CC$-injective)
and $Y \in \CC$.
If $\CC$ has enough projectives and enough injectives, and
if these coincide (i.e. an object is $\CC$-projective if and only if
it is $\CC$-injective), then
$\CC$ is called a {\it Frobenius category}.

For such a Frobenius category $\CC$ we can define its 
{\it stable category} ${\stCC}$.  
The objects of ${\stCC}$ are the 
same as the objects of $\CC$, and the morphism spaces 
$\Hom_{\stCC}(X,Y)$ are the morphism
spaces in $\CC$ modulo morphisms factoring through 
$\CC$-projective-injective objects.
The category ${\stCC}$ is a triangulated 
category in a natural way \cite{H2}, where
the shift is given by the relative inverse syzygy functor
$
\Omega^{-1}\colon {\stCC} \to {\stCC}.
$
(Recall that the syzygy functor $\Omega$ maps an object 
to the kernel of its projective cover. This induces an 
auto-equivalence of the stable category.)
For all $X$ and $Y$ in $\CC$ there
is a functorial isomorphism
\begin{equation}\label{functorial1}
\Ext_\CC^1(X,Y) \cong \Hom_{\stCC}(X,\Omega^{-1}(Y)).
\end{equation}
The category $\CC$ is called \emph{stably $2$-Calabi-Yau} if the stable category
$\stCC$ is a {\it $2$-Calabi-Yau category}, that is, if
for all $X,Y \in \CC$ there is
a functorial isomorphism
\begin{equation}\label{functorial2}
\Ext_\CC^1(X,Y) \cong {\rm D}\Ext_\CC^1(Y,X),
\end{equation}
where ${\rm D}$ denotes the duality for vector spaces.
We can now state some important properties of the categories $\CC_w$.
\begin{theorem}[\cite{BIRS,GLSUni1}]\label{main1}
The categories $\CC_w$ are Frobenius categories. They are stably $2$-Calabi-Yau. 
The indecomposable $\CC_w$-projective-injective modules,
are the indecomposable direct summands of $I_w$. 
\end{theorem}

\subsection{Cluster-tilting modules}
\label{sect4.5}

Let $\CC$ be a subcategory of $\nil(\L)$ closed under extensions, direct sums and direct summands. 
For an object $T$ of $\CC$
we denote by $\add(T)$ the \emph{additive envelope} of $T$, that is,
the full subcategory whose objects are finite direct sums
of direct summands of $T$.
We say that
\begin{itemize}
\item
$T$ is $\CC$-{\it maximal rigid} if 
$
\Ext_\L^1(T \oplus X,X) = 0
$
with $X \in \CC$ implies $X \in \add(T)$;
\item
$T$ is a $\CC$-{\it cluster-tilting module} if 
$
\Ext_\L^1(T,X) = 0
$
with $X \in \CC$
implies $X \in \add(T)$.
\end{itemize}
Since $\Ext_\L^1(T \oplus X,X) \cong \Ext_\L^1(T,X) \oplus \Ext_\L^1(X,X)$,
we see that the second property is a priori stronger than the first.

\begin{theorem}[\cite{BIRS,GLSUni1}]
\label{main9}
For a rigid $\L$-module $T$ in $\CC_w$ the following are equivalent:
\begin{itemize}

\item[(i)]
$T$ has $r$ pairwise non-isomorphic indecomposable direct summands;

\item[(ii)]
$T$ is $\CC_w$-maximal rigid;

\item[(iii)]
$T$ is a $\CC_w$-cluster-tilting module.

\end{itemize}
\end{theorem}

Let $T = T_1\oplus \cdots \oplus T_r$ be a $\CC_w$-cluster-tilting module, with each
summand $T_i$ indecomposable. 
Consider the endomorphism algebra $A_T:=\End_\L(T)^{\rm op}$. This is a basic algebra,
with indecomposable projective modules 
\[
P_{T_i} := \Hom_\L(T,T_i),\qquad (1\le i\le r).
\]
The simple $A_T$-modules are the heads $S_{T_i}$ of the projectives $P_{T_i}$.
One defines a quiver $\G_T$ with vertex set $\{1,\ldots,r\}$, and 
$d_{ij}$ arrows from $i$ to $j$, where 
\[
d_{ij} := \dim \Ext^1_{A_T}(S_{T_i},S_{T_j}).
\] 
This is known as the \emph{Gabriel quiver} of $A_T$. 
%

Let now $\i=(i_r,\ldots,i_1)$ be a reduced expression for $w$.
Following \cite{BFZ}, we define a quiver
$\G_\i$ as follows.
The vertex set of $\G_\i$ is equal to $\{1,\ldots,r\}$.
For $1 \le k \le r$, let 
\begin{align*}\label{kminus}
k^- &:= \max\left(\{0\}\cup \{1 \le s \le k-1 \mid i_s = i_k \}\right),\\
k^+ &:= \min\left(\{ k+1 \le s \le r \mid i_s = i_k \} \cup \{r+1\}\right).
\end{align*}
For $1 \le s,\,t \le r$ such that $i_s \not = i_t$, 
there are $|c_{i_s,i_t}|$ arrows from $s$ to $t$ provided
$t^+ \ge s^+ > t > s$.
(Here, as in \S\ref{sect3.1}, the $c_{ij}$'s are the entries of the Cartan matrix.)
These are called the {\it ordinary arrows} of $\G_\i$.
Furthermore, for each $1 \le s \le r$ there is an arrow $s \to s^-$
provided $s^- > 0$.
These are the {\it horizontal arrows} of $\G_\i$.

\begin{theorem}[\cite{BIRS,GLSUni1}]
\label{main6}
The module $V_\i$ is a $\CC_w$-cluster-tilting module, 
and we have $\G_{V_\i} = \G_\i$.
\end{theorem}

\begin{example}\label{Example4.8}
We continue Example~\ref{Example4.4}. For $\i=(2,1,2,1)$, the Gabriel quiver
$\G_\i$ of $\End_\L(V_\i)^{\rm op}$ is:
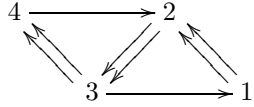
\[
\xymatrix@-1.5pc{
4 \ar@{->}[rrrr] \ar@<.5ex>@{<-}[rrdd]\ar@<-.5ex>@{<-}[rrdd]& & & & 
\ar@<.5ex>@{->}[lldd]\ar@<-.5ex>@{->}[lldd]2 \ar@<.5ex>@{<-}[rrdd]\ar@<-.5ex>@{<-}[rrdd]\\
&&&\\
&& 3 \ar@{->}[rrrr] & & & & 1 
}
\]
Note that the $\CC_w$-projective summands of $V_\i$ correspond to the leftmost
vertices of each row of $\G_\i$.
\end{example}

\subsection{Mutations of cluster-tilting modules}

Let $T=T_1\oplus \cdots \oplus T_r$ be a $\CC_w$-cluster-tilting module,
and let $\G_T$ be the corresponding quiver defined in \S\ref{sect4.5}.
Define 
\[
b_{ij} := (\text{number of arrows } j\to i \text{ in } \G_T)
- (\text{number of arrows } i\to j \text{ in } \G_T). 
\]
Clearly each indecomposable $\CC_w$-projective-injective module is
a direct summand of $T$. In the sequel we will not need the arrows 
of $\G_T$ between the vertices corresponding to these projective-injective summands.
\begin{theorem}[\cite{BIRS,GLSUni1}]\label{th_mutation}
Let $T_k$ be a non-projective indecomposable direct summand of $T$. 
\begin{itemize}
\item[(i)] There exists a unique indecomposable module $T_k^*\not\cong T_k$
such that $(T/T_k)\oplus T_k^*$ is a $\CC_w$-cluster-tilting module.
We call $(T/T_k)\oplus T_k^*$ the \emph{mutation of $T$ in direction $k$},
and denote it by $\mu_k(T)$. 
\item[(ii)] The quiver $\G_{\mu_k(T)}$ is equal to the  
mutation $\mu_k$ of $\G_T$ in the sense of Fomin and Zelevinsky 
(ignoring the arrows between projective-injective vertices). 
\item[(iii)] We have $\dim\Ext^1_\L(T_k,T_k^*) = 1$.
Let 
\[
0\to T_k \to T'_k \to T_k^* \to 0 \quad\text{and}\quad 
0\to T_k^* \to T''_k \to T_k \to 0 
\]
be non-split short exact sequences. Then 
\[
 T'_k \cong \bigoplus_{b_{jk}<0} T_j^{-b_{jk}},\qquad
 T''_k \cong \bigoplus_{b_{jk}>0} T_j^{b_{jk}}.
\]
\end{itemize}
\end{theorem}

\begin{remark}
Let $\underline{T}$ be any cluster-tilting object of the stable category 
$\stCC_w$, coming from a cluster-tilting module $T$ in the mutation class of a module
$V_\i$. It was shown by Buan, Iyama, Reiten and Smith \cite{BIRS2} that the endomorphism
algebra of $\underline{T}$ is the Jacobian algebra of a quiver with potential.
Jacobian algebras and their mutations have been introduced
by Derksen, Weyman and Zelevinsky \cite{DWZ1}, and they have been used 
in \cite{DWZ2} to prove several important conjectures on cluster algebras by 
Fomin and Zelevinsky \cite{FZ4}.
\end{remark}

\subsection{A distinguished sequence of mutations}\label{section4.7}

For $1\le k\le l\le r$ such that $i_k=i_l=i$, we have a natural embedding 
of $\L$-modules $V_{k^-} \to V_{l}$.
Following \cite[\S9.8]{GLSKM}, we define $M[l,k]$ as the cokernel 
of this embedding, that is,
\[
 M[l,k] := V_{l} / V_{k^-}.
\]
In particular, we set $M_k:= M[k,k]$, and
\[
M_\i := M_1\oplus \cdots \oplus M_r. 
\] 
We will use the convention that $M[l,k]=0$ if $k>l$.
Every module $M[l,k]$ is indecomposable and rigid.
But note that $M_\i$ is \emph{not} rigid.
Define 
\begin{align*}\label{kmin}
k_{\min} &:= \min\{1 \le s \le r \mid i_s = i_k \},\\
k_{\max} &:= \max\{ 1 \le s \le r \mid i_s = i_k \}.
\end{align*}
Then $V_k = M[k,k_{\min}]$ corresponds to an \emph{initial interval}.
The direct sum of all modules $M[k_{\max},k]$ corresponding to 
\emph{final intervals} is also a $\CC_w$-cluster-tilting module, denoted 
by $T_\i$. We number the summands of $T_\i$ as follows: 
$$
T_k := 
\begin{cases}
V_k & \text{if $k^+ = r+1$},\\
M[k_{\rm max},k^+] & \text{otherwise}.
\end{cases}
$$
This numbering ensures that, for a non-projective $V_k$ we have $\Omega_w^{-1}(V_k) = T_k$,
where $\Omega_w$ denotes the relative syzygy functor of $\stCC_w$. 
\begin{example}\label{Example4.10}
We continue Example~\ref{Example4.4} and Example~\ref{Example4.8}.
In this case we have $M_1 = V_1$, $M_2=V_2$, and the modules $M_3$
and $M_4$ were already introduced in Example~\ref{Example4.4}. 
The indecomposable direct summands of $T_\i$ are then 
\[
T_1 = M_3,\quad
T_2 = M_4,\quad
T_3 = V_3,\quad
T_4 = V_4, 
\]
and the quiver $\G_{T_\i}$ is 
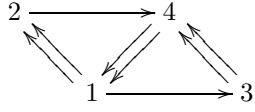
\[
\xymatrix@-1.5pc{
2 \ar@{->}[rrrr] \ar@<.5ex>@{<-}[rrdd]\ar@<-.5ex>@{<-}[rrdd]& & & & 
\ar@<.5ex>@{->}[lldd]\ar@<-.5ex>@{->}[lldd]4 \ar@<.5ex>@{<-}[rrdd]\ar@<-.5ex>@{<-}[rrdd]\\
&&&\\
&& 1 \ar@{->}[rrrr] & & & & 3 
}
\]
Note that this quiver is the same as $\G_{V_\i}$, but the vertices corresponding 
to the $\CC_w$-projective summands are now at the right end of each row. 
\end{example}

It was shown in \cite[\S13.1]{GLSKM} that there is an explicit sequence
of mutations of $\CC_w$-cluster-tilting modules starting from $V_\i$ and
ending in $T_\i$. 
This sequence of mutations plays an important role in several of our constructions,
so we want to describe it in some detail.
For $j\in I$ and $1 \le k \le r+1$, we set
\[
k[j] := |\{ 1 \le s \le k-1 \mid i_s = j \}|,
\qquad
t_j := (r+1)[j].
\]
Thus $t_j$ is the number of occurences of $j$ in $\i$.
Our sequence consists of 
\[
\sum_{j\in I} \frac{t_j(t_j-1)}{2}
\] 
mutations,
which we conveniently group into $r$ steps.
We start from $V_\i$ and its quiver $\G_\i$.
Note that $\G_\i$ is naturally displayed on $n$ rows,
where all vertices $k$ such that $i_k=j$ are sitting
on row $j$, with their labels increasing from right to left
(see Example~\ref{Example4.8}).
At Step $1$, we perform $t_{i_1}-1-1[i_1]=t_{i_1}-1$ mutations
at consecutive vertices sitting on row $i_1$ of $\G_\i$, starting 
from the rightmost vertex. Next, at Step~$k=2,\ldots,r$, 
we perform $t_{i_k}-1-k[i_k]$ mutations
at consecutive vertices sitting on row $i_k$, starting 
from the rightmost vertex. 
Now we claim:
\begin{theorem}[{\cite[\S13]{GLSKM}}]\label{Theorem4.11}
\begin{itemize}
\item[(i)] The above sequence of mutations applied to $V_\i$ gives $T_\i$.
\item[(ii)] All the indecomposable direct summands of the $\CC_w$-cluster-tilting
modules occuring in this sequence are of the form $M[l,k]$.
\item[(iii)] Each step of this sequence consists of the mutation
of a module $M[d^-,b]$ into a module $M[d,b^+]$, for some 
$1\le b < d \le r$ with $i_b = i_d = i$.
The corresponding pair of short exact sequences is 
\[
0 \to M[d^-,b] \to M[d^-,b^+] \oplus M[d,b] \to M[d,b^+] \to 0, 
\]
\[
0 \to M[d,b^+] \to \bigoplus_{j\not = i} M[d^-(j),b^+(j)]^{\oplus |c_{ij}|} \to M[d^-,b] \to 0,
\]
where for $1\le k\le r$, we set 
\[
k^-(j):=\max\{0, s < k \mid i_s=j\},
\quad k^+(j):=\min\{r+1, k< s  \mid i_s=j\}.
\]
\item[(iv)] Every module $M[l,k]$ with $1\le k\le l\le r$ and $i_k=i_l$
arises in this sequence.
\end{itemize}
\end{theorem}

\begin{example} \label{Example4.12}
We continue Example~\ref{Example4.10}.
We have $t_1=t_2=2$, hence the sequence consists of only $1+1=2$ mutations.

{\it Step 1}: we perform mutation $\mu_1$. The two short exact sequences are
\[
0\to V_1 \to V_3 \to T_1 \to 0,\qquad
0\to T_1 \to V_2^{\oplus 2} \to V_1 \to 0. 
\]
Note that we have 
\[
V_1 = M[3^-,1],\quad T_1 = M[3, 1^+],\quad V_3 = M[3,1],\quad V_2 = M[3^-(2),1^+(2)].
\]

{\it Step 2}: we perform mutation $\mu_2$. The two short exact sequences are
\[
0\to V_2 \to V_4 \to T_2 \to 0,\qquad
0\to T_2 \to T_1^{\oplus 2} \to V_2 \to 0. 
\]
Note that we have 
\[
V_2 = M[4^-,2],\quad T_2 = M[4, 2^+],\quad V_4 = M[4,2],\quad T_1 = M[4^-(1),2^+(1)].
\]
For more complicated examples, we refer to \cite[\S13]{GLSKM}.
\end{example}

We conclude this section by noting that
the functions $\varphi_{M[l,k]}$ associated with the 
modules $M[l,k]$ are restrictions to $N_+$ of some generalized minors in the sense of \S\ref{sect3.3}.
In particular in type $A_n$, they are nothing but ordinary minors 
of a uni\-triangular matrix of size $n+1$.
Using Theorem~\ref{Th_mult}, we can convert the above mutation sequence 
into a sequence of determinantal identities, and thus recover certain
identities of Fomin and Zelevinsky \cite{FZTP}.

\section{Cluster structures on coordinate rings}
\label{sect5}

We can now use the categories $\CC_w$ of \S\ref{sect4.3} to produce 
some cluster algebras, and show that they are isomorphic 
to the coordinate rings of $N(w)$ and $N^w$. 
Our construction readily implies that the cluster monomials 
are contained in the dual semicanonical basis $\SC^*$ of $\C[N]$. 

\subsection{From categories to cluster algebras}

We say that a $\CC_w$-cluster-tilting module is \emph{reachable} if it can 
be obtained from $V_\i$ as the result of a (finite) sequence of mutations.
It is known \cite[Lemma II.4.2]{BIRS} that if $\j$ is another reduced word for~$w$, then $V_\j$
is reachable from $V_\i$, hence this notion does not depend on the choice of $\i$.
It is an open problem whether every $\CC_w$-cluster-tilting module is reachable.
More generally, we call a $\L$-module \emph{reachable} if it is isomorphic to
a direct summand of a reachable $\CC_w$-cluster-tilting module.

Let $\RR(\CC_w)$ be the subalgebra of $\C[N]$ generated by the 
$\varphi_{T_k}\ (1\le k\le r)$
where $T=T_1\oplus\cdots\oplus T_r$ runs over all reachable $\CC_w$-cluster-tilting modules.
Let $\A_\i$ denote the cluster algebra defined by the initial seed
$((y_1,\ldots,y_r),\G_\i)$, in which the variables $y_k$ corresponding
to $\CC_w$-projective-injective vertices are frozen.

\begin{theorem}[\cite{GLSKM}]\label{Theorem5.1}
\begin{itemize}
\item[(i)]
There is a unique algebra isomorphism $\iota$ from $\A_\i$ to
$\RR(\CC_w)$ such that 
\[
\iota(y_k)=\varphi_{V_k}, \qquad (1\le k\le r).
\]
\item[(ii)]
If we identify the two algebras $\A_\i$ and
$\RR(\CC_w)$ via $\iota$, then the clusters of
$\A_\i$ are identified with the $r$-tuples
$(\varphi_{T_1},\ldots,\varphi_{T_r})$, where $T$
runs over all reachable $\CC_w$-cluster-tilting modules.
Moreover, all cluster monomials belong to 
the dual semicanonical basis 
$\SC^*$ of $\C[N]$.
\end{itemize} 
\end{theorem}
Theorem~\ref{Theorem5.1} gives a Lie-theoretic realization
of a large class of cluster algebras. Its proof relies on Theorem~\ref{Th_mult}
and Theorem~\ref{th_mutation}. 
As an application, we can compute
the Euler characteristics $\chi(\F_{X,\i})$
for all modules $X$ in the additive closure of a reachable
$\CC_w$-cluster-tilting module,
and all composition series types $\i$, using the Fomin-Zelevinsky
mutation formula (\ref{FZmutation_form}). 
Equivalently, we have an algorithm for computing the elements
$\varphi_X$ of $\SC^*$ corresponding to any reachable rigid module.

The next theorem describes $\C$-bases of the above cluster algebras,
and shows that they are polynomial rings. 
Recall from \S\ref{section4.7} the $\L$-module 
\[
M_\i = M_1\oplus\cdots\oplus M_r\in \CC_w.
\]

\begin{theorem}[\cite{GLSKM}]\label{basesthm}
\begin{itemize}

\item[(i)]
The cluster algebra $\RR(\CC_w)$ is a polynomial ring in
$r$ variables.
More precisely, we have
$$
\RR(\CC_w) = \C[\varphi_{M_1},\ldots,\varphi_{M_r}]. 
$$
\item[(ii)]
The set 
$
\mathcal{P}_\i^*:=\left\{ \varphi_M \mid M \in \add(M_\i) \right\}
$ 
is a $\C$-basis of $\RR(\CC_w)$.
\item[(iii)]
The subset 
$\SC_w^* := \SC^* \cap \RR(\CC_w)$
of the dual semicanonical basis
of $\C[N]$
is a $\C$-basis of $\RR(\CC_w)$ containing all cluster
monomials.
\end{itemize}
\end{theorem}

The bases given by Theorem~\ref{basesthm} (ii) and (iii)
are called \emph{dual PBW basis} and \emph{dual semicanonical basis}
of $\RR(\CC_w)$, respectively.
The proof of this theorem uses Theorem~\ref{Theorem4.11} for
showing that $\C[\varphi_{M_1},\ldots,\varphi_{M_r}] \subseteq \RR(\CC_w)$.
The reverse inclusion is obtained by proving that $\mathcal{P}_\i^*$
is a subset of a basis of $\C[N]$ dual to a Poincar\'e-Birkhoff-Witt
basis of $U(\n)$ (see \cite[\S15]{GLSKM}).

\subsection{Coordinate rings of unipotent subgroups and unipotent cells}

Using the fact that $\mathcal{P}_\i^*$ is a subset of an appropriate
dual PBW-basis, one can show that the functions $\varphi_{M_k}$
are $N'(w)$-invariant. It then follows from Proposition~\ref{Prop3.2}
that we can relate the subalgebra $\RR(\CC_w)$ of $\C[N]$ to 
the coordinate ring of $N(w)$, as explained in the next theorem.

\begin{theorem}[\cite{GLSKM}]\label{thm5.3}
The cluster algebra $\RR(\CC_w)$ coincides with the invariant subring 
$\C[N]^{N'(w)}$, and is naturally isomorphic to $\C[N(w)]$.  
\end{theorem}
 
As a result, we have obtained that the coordinate ring $\C[N(w)]$
has the structure of a cluster algebra. Moreover, its clusters are
in one-to-one correspondence with reachable $\CC_w$-cluster-tilting objects.
The frozen cluster variables are the $\varphi_{I_i}$, where $I_i$ 
is the indecomposable $\CC_w$-projective-injective module with simple socle~$S_i$.
It turns out that
\[
 \varphi_{I_i} = D_{\varpi_i,\,w^{-1}(\varpi_i)}, \qquad (i\in I), 
\]
where we denote by $D_{u(\varpi_i),v(\varpi_i)}$ the restriction to
$N_+$ of the generalized minor $\De_{u(\varpi_i),v(\varpi_i)}$.
The coefficient ring of the cluster algebra $\C[N(w)]$ is therefore
the polynomial ring in the frozen variables $D_{\varpi_i,w^{-1}(\varpi_i)}$.

By Proposition~\ref{Prop3.5}, these frozen variables are precisely
the functions that need to be inverted in order to pass from
$\C[N(w)]$ to the coordinate ring $\C[N^w]$ of the unipotent cell
$N^w$.
Therefore $\C[N^w]$ has almost the same cluster algebra structure
as $\C[N(w)]$. 
The only difference 
is in the coefficient ring, which in the case of $N^w$ should be 
taken as the ring of \emph{Laurent polynomials} in the generalized minors
$D_{\varpi_i,\,w^{-1}(\varpi_i)}$. Hence we have

\begin{theorem}[\cite{GLSKM}]\label{thm5.4}
The cluster algebra $\widetilde{\RR}(\CC_w)$ obtained from 
$\RR(\CC_w)$ by formally inverting all the frozen variables
is naturally isomorphic to $\C[N^w]$. 
\end{theorem}
When $\g$ is finite-dimensional, the cluster algebra structure
on the unipotent cell $N^w$ given by Theorem~\ref{thm5.4} was
already obtained by Berenstein, Fomin and Zelevinsky \cite{BFZ},
by a completely different method.
In fact they treated the more general case of \emph{double Bruhat cells}
\[
G^{v,w}:=(BvB) \cap (B_-wB_-),\qquad (v,w \in W). 
\]
These varieties $G^{v,w}$ carry a free action of the torus $H$
by left (or right) multiplication, and quotienting this action
gives the \emph{reduced double Bruhat cells}
\[
L^{v,w}:=(NvN) \cap (B_-wB_-),\qquad (v,w \in W). 
\]
One can then show that going from $G^{e,w}$ to
$L^{e,w}=N^w$ 
only modifies the coefficient ring of the cluster algebra. 
However, our approach shows that $\C[N^w]$ is a genuine cluster
algebra (not only an \emph{upper} cluster algebra in the sense of \cite{BFZ}),
and it shows that cluster monomials belong to the dual semicanonical basis.

In the general Kac-Moody case, Theorem~\ref{thm5.4} was conjectured 
by Buan, Iyama, Reiten and Scott \cite[Conjecture III.3.1]{BIRS}.  
Theorem~\ref{thm5.4} has also been extended to the non simply-laced
case and adaptable Weyl group element $w$
by Demonet \cite{D}.

\begin{example}
We continue Example~\ref{Example3.7}. 
The full subquiver obtained by erasing vertices $3$ and $4$ of the 
quiver of Example~\ref{Example4.8} is the Kronecker quiver with
two arrows. 
It follows that $\C[N(w)]$ is a rank 2 acyclic cluster algebra 
of affine type $\widetilde{A}_1$,
that is, a version of the cluster algebra $\A_{\Sigma}$ of Example~\ref{ex2.1},
with $a=2$ and two additional frozen variables.
It has infinitely many clusters and cluster variables.
 
Using the notation of Example~\ref{Example4.4} and Example~\ref{Example4.10},
the dual PBW generators $\varphi_{M_i}$ evaluated at an element $x\in N(w)$ are expressed, in terms of 
the coordinate functions $b_k$ on $\C[N(w)]$
introduced in Example~\ref{Example3.7}, as
\begin{equation}\label{eq6}
\varphi_{M_1} = b_0,\quad \varphi_{M_2} = -b_1,\quad  \varphi_{M_3} = b_2, \quad
\varphi_{M_4} = -b_3.
\end{equation}
Our initial cluster for
$\C[N(w)]$ is $(\varphi_{V_1}, \varphi_{V_2}, \varphi_{V_3}, \varphi_{V_4})$,
where
\begin{equation}\label{eq7}
\varphi_{V_1} = b_0,\quad \varphi_{V_2} = -b_1,\quad  \varphi_{V_3} = b_0b_2-b_1^2, \quad
\varphi_{V_4} = b_1b_3-b_2^2.
\end{equation}
In agreement with Example~\ref{Example4.12},
one can check the exchange relations corresponding to a mutation at vertex $1$
followed by a mutation at vertex 2:
\[
\varphi_{V_1}\varphi_{M_3} = \varphi_{V_3}+\varphi_{V_2}^2,\qquad
\varphi_{V_2}\varphi_{M_4} = \varphi_{V_4}+\varphi_{M_3}^2, 
\]
If instead we start mutating at vertex $2$, we get the new cluster variable
\[
\frac{b_0^2(b_1b_3-b_2^2)+(b_0b_2-b_1)^2}{(-b_1)} = 2b_0b_1b_2-b_1^3-b_0^2b_3. 
\]

We may also evaluate $\varphi_{V_i}$ and $\varphi_{M_i}$ at an
arbitrary point 
$x = \begin{pmatrix}
a & b \\
c & d 
\end{pmatrix}$ of~$N_+$. We obtain the following determinantal expressions
in terms of the coordinate functions $b_k$ and $d_k$ of Example~\ref{Example3.6}:
\[
\begin{array}{ll}
\varphi_{M_1}=\varphi_{V_1} = b_0,
&\varphi_{M_2}=\varphi_{V_2} = 
\left|\begin{matrix}
b_0 & b_1 \\
1 & d_1 
\end{matrix}
\right|,
\\[5mm]
\varphi_{V_3} = 
\left|\begin{matrix}
b_0 & b_1 &b_2\\
1 & d_1 & d_2 \\
0 & b_0 & b_1
\end{matrix}
\right|, 
&
\varphi_{M_3} = 
\left|\begin{matrix}
b_0 & b_1 &b_2\\
1 & d_1 & d_2 \\
0 & 1 & d_1
\end{matrix}
\right|,
\\[8mm]
\varphi_{V_4} = 
\left|\begin{matrix}
b_0 & b_1 &b_2 &b_3\\
1 & d_1 & d_2 &d_3\\
0 & b_0 & b_1 &b_2\\
0&1&d_1&d_2
\end{matrix}
\right|, 
&
\varphi_{M_4} = 
\left|\begin{matrix}
b_0 & b_1 &b_2&b_3\\
1 & d_1 & d_2 &d_3\\
0 & 1 & d_1 &d_2\\
0 & 0 & 1 &d_1
\end{matrix}
\right|.
\end{array}
\]
(These formulas give back equations (\ref{eq6}), (\ref{eq7}), if we specialize
$d_1 = d_2 = d_3 = 0$.)
Using the description of $N'(w)$ from Example~\ref{Example3.7},
it is easy to check that these determinantal expressions are $N'(w)$-invariant,
as claimed in Theorem~\ref{thm5.3}.

For more general determinantal and combinatorial evaluations of functions 
$\varphi_X$ in the case of $\slchap_2$, see \cite{Sc}. 

\end{example}

\section{Chamber Ansatz}\label{sectCA}

We now come back to the problem of computing the factorization
parameters $t_k$ of a point $x_\i(\t)$ of the unipotent cell $N^w$
(see \S\ref{sect3.5}). To do so, we introduce a new distinguished
cluster-tilting object of $\CC_w$.
Let $\Omega_w$ be the syzygy functor of the stable category $\stCC_w$.
Thus, for a module $X\in\CC_w$, we have a short exact sequence
\[
0\to \Omega_w(X) \to P(X) \to X \to 0, 
\]
where $P(X)$ denotes the projective cover of $X$ in $\CC_w$.
Note that for a $\CC_w$-projective-injective object $X$, $\Omega_w(X)=0$.
Define
\[
W_\i := \Omega_w(V_\i) \oplus I_w. 
\]
Since $\Omega_w$ is an auto-equivalence of $\stCC_w$, we see that
$W_\i$ is indeed cluster-tilting.

Recall from \S\ref{section4.7} that we have another $\CC_w$-cluster-tilting
module $T_\i$, related to $V_\i$ by a distinguished sequence of mutations.
It is easy to check that 
\[
V_\i = \Omega_w(T_\i) \oplus I_w. 
\]
This implies that, given a reachable $\CC_w$-cluster-tilting module $M$,
we have an ``algorithm'' for producing a sequence of mutations from
$M$ to $\Omega_w(M) \oplus I_w$ (see \cite[Proposition 13.4]{GLSKM}).
In particular, $\Omega_w(M) \oplus I_w$ is again reachable.
Hence $W_\i$ is reachable, and the functions $\varphi_{W_k}$
associated with its indecomposable direct summands $W_k\ (1\le k\le r)$
form a new cluster of the cluster algebra $\C[N^w]$. 
For $k=1,\ldots,r$, define
\begin{equation}
\varphi'_{V_{k}} := \frac{\varphi_{W_{k}}}{\varphi_{P(V_{k})}},
\end{equation}
a Laurent monomial in the $\varphi_{W_{k}}$ (since $\add(W_\i)$ contains all $\CC_w$-projective-injective modules), 
and put
\begin{equation}\label{ch5}
C_{k} := \frac{1}{\vph_{V_{k}}'\vph_{V_{k^-(i_k)}}'} \cdot
\prod_{j\not = i_k} \left(\vph_{V_{k^-(j)}}'\right)^{|c_{i_k,j}|}.
\end{equation}
Here we recall that $k^-(j) := \max\{ 0,1 \le s \le k-1 \mid i_s = j \}$ and
$V_{0}$ is by convention the zero module.
We can now state the Chamber Ansatz formula, 
which expresses the rational function $t_k$ as an explicit Laurent monomial
in the cluster variables~$\varphi_{W_k}$.

\begin{theorem}\label{THM5}
For $1 \le k \le r$ and $\mathbf{x}_\i(\t) = x_{i_r}(t_r)\cdots x_{i_1}(t_1)$ we have 
\[
t_k = C_{k}(\mathbf{x}_\i(\t)).
\]
\end{theorem}

\begin{example}\label{example6.2}
We continue Example~\ref{Example3.6}.
First we compute $W_\i$. 
We have short exact sequences
$$
0 \to W_1\to P(V_1) \to V_1 \to 0
\qquad\text{ and }\qquad
0 \to W_2 \to P(V_2) \to V_2 \to 0.
$$
where $P(V_1) \cong V_3^{\oplus 3}$ and $P(V_2) \cong V_3^{\oplus 2}$.
Thus, using the same convention as in Example~\ref{Example4.4}, we can
represent $W_1$ and $W_2$ as
\[
\begin{array}{l}
W_1\ =\ \vcenter{\def\objectstyle{\scriptstyle}\def\labelstyle{\scriptstyle}
\xymatrix@-1.5pc{
1\ar@{-}[rd]&&\ar@{-}[ld]1\ar@{-}[rd]&&\ar@{-}[ld]1\ar@{-}[rd]&&\ar@{-}[ld]1&&1\ar@{-}[rd]&&\ar@{-}[ld]1\ar@{-}[rd]&&\ar@{-}[ld]1\ar@{-}[rd]&&\ar@{-}[ld]1\\
 &2\ar@{.}[rd]& &\ar@{.}[ld]2\ar@{.}[rrrrd]& &2\ar@{.}[rrd]& & & &\ar@{.}[lld]2& &\ar@{.}[lllld]2\ar@{.}[rd]& &\ar@{.}[ld]2\\
 & &1& & & & &1& & & & &1\\ 
}}\\[10mm]
W_2\ =\ \vcenter{\def\objectstyle{\scriptstyle}\def\labelstyle{\scriptstyle}
\xymatrix@-1.5pc{
1\ar@{-}[rd]& &\ar@{-}[ld]1\ar@{-}[rd]& &\ar@{-}[ld]1\ar@{-}[rd]& &\ar@{-}[ld]1\\  
 &2\ar@{.}[rd]& &\ar@{.}[ld]2\ar@{.}[rd]& &\ar@{.}[ld]2& \\
 & &1& &1& & \\
}}
\end{array}
\]
Putting $\sx_\i(\t) = x_2(t_4)x_1(t_3)x_2(t_2)x_1(t_1)$, one can then calculate
\[
\begin{array}{lcl}
\varphi_{V_1}(\sx_\i(\t))   &=& t_3 + t_1,\\[2mm]
\varphi_{V_2}(\sx_\i(\t))   &=& t_4(t_3^2+ 2t_3t_1+t_1^2)+ t_2 t_1^2,\\[2mm]
\varphi_{V_3}(\sx_\i(\t))   &=& t_3 t_2^2t_1^3,\\[2mm]
\varphi_{V_4}(\sx_\i(\t))   &=& t_4 t_3^2t_2^3t_1^4,\\[2mm]
\varphi_{W_1}(\sx_\i(\t))   &=& t_3^3t_2^6t_1^8,\\[2mm]
\varphi_{W_2}(\sx_\i(\t))   &=& t_3^2t_2^3t_1^4.
\end{array} 
\]
Noting that 
\[
\varphi_{W_1} = \varphi'_{V_1}\varphi_{V_3}^3,
\quad
\varphi_{W_2} = \varphi'_{V_2}\varphi_{V_3}^2,
\quad
\varphi_{W_3} = \varphi_{V_3} = \frac{1}{\varphi'_{V_3}},
\quad
\varphi_{W_4} = \varphi_{V_4} = \frac{1}{\varphi'_{V_4}},
\]
we thus get, in agreement with Theorem~\ref{THM5},
\[
\begin{array}{lclcl}
t_1&=&\displaystyle\frac{1}{\varphi'_{V_1}}(\sx_\i(\t))&
=&\displaystyle\frac{\varphi_{W_3}^3}{\varphi_{W_1}}(\sx_\i(\t)),\\[5mm]
t_2&=&\displaystyle\frac{(\varphi'_{V_1})^2}{\varphi'_{V_2}}(\sx_\i(\t))&
=&\displaystyle\frac{\varphi_{W_1}^2}{\varphi_{W_2}\varphi_{W_3}^4}(\sx_\i(\t)),\\[5mm]
t_3&=&\displaystyle\frac{(\varphi'_{V_2})^2}{\varphi'_{V_3}\varphi'_{V_1}}(\sx_\i(\t))&
=&\displaystyle\frac{\varphi_{W_2}^2}{\varphi_{W_1}}(\sx_\i(\t)),\\[5mm]
t_4&=&\displaystyle\frac{(\varphi'_{V_3})^2}{\varphi'_{V_4}\varphi'_{V_2}}(\sx_\i(\t))&
=&\displaystyle\frac{\varphi_{W_4}}{\varphi_{W_2}}(\sx_\i(\t)). 
\end{array}
\]
As already mentioned in Example~\ref{Example3.6}, the numerators and denominators of
these rational functions differ from those of (\ref{eqnotCA}), which are \emph{not} cluster
monomials, and which are \emph{not} $N'(w)$-invariant.
For example the denominator $a_1$ of (\ref{eqnotCA}) is the restriction to 
$N^w$ of the function $\varphi_X$ of $\C[N]$, where $X$ is the $2$-dimensional
$\L$-module described at the end of Example~\ref{Example4.4}, which is not
an object of $\CC_w$.
\end{example}

\section{Quantum cluster structures on quantum coordinate rings}

The coordinate ring $\C[N(w)]$ has a quantum analogue $U_q(\n(w))$
introduced by De Concini, Kac and Procesi. On the other hand, Berenstein and
Zelevinsky \cite{BZ} have introduced quantum analogues of cluster
algebras. In this section we explain that $U_q(\n(w))$ has a quantum
cluster algebra structure obtained by $q$-deforming in an appropriate 
way the cluster algebra structure of $\C[N(w)]$.

\subsection{Quantum cluster algebras}

The guiding principle is to replace the Laurent polynomial rings generated by cluster
variables of any given cluster by  \emph{quantum tori}, that is, to require that the 
corresponding quantum cluster variables are pairwise $q$-commutative. 
A certain compatibility condition then ensures that all these quantum tori 
can be glued together to form a flat deformation of the original cluster algebra.
It is important to note that a given cluster algebra may have several non-isomorphic
$q$-deformations, and it may also have no $q$-deformation at all (if its exchange
matrix does not have maximal rank).

Instead of repeating the  definition of a quantum cluster algebra
(for which we refer to \cite{BZ}),
let us describe an example constructed from the Frobenius category $\CC_w$ of
\S\ref{sect4.3}. 
For $M, N \in \CC_w$ let us write for short
$[M,N]:= \dim\Hom_\L(M,N)$. 
Recall the $\CC_w$-cluster-tilting object $V_\i=V_1\oplus\cdots\oplus V_r$, and define
\[
\la_{ij} := [V_i, V_j] - [V_j, V_i], \qquad (1\le i,j\le r).
\]
Let $q$ denote an indeterminate over $\Q$. We introduce the quantum torus
\[
\T := \Q(q)\langle Y_1^{\pm 1},\ldots , Y_r^{\pm 1}\rangle
\] 
whose generators $Y_k$ obey the $q$-commutation relations
\[
Y_iY_j = q^{\la_{ij}} Y_jY_i.
\]
For $R = V_1^{a_1}\oplus \cdots \oplus V_r^{a_r}\in \add(V_\i)$, set 
$Y_R := q^{-\a(R)} Y_1^{a_1}\cdots Y_r^{a_r}$, where
\[
\a(R) := \sum_{1\le i<j\le r}a_ia_j[V_i,V_j] + \sum_{1\le i\le r} \frac{a_i(a_i-1)}{2} [V_i,V_i]. 
\]
In particular, $Y_{V_i}=Y_i$.
It is easy to check that for any modules $R, S \in \add(V_\i)$, 
we have 
\[
Y_R Y_S = q^{[R,S]} Y_{R\oplus S}.
\]
\begin{theorem}[\cite{GLSquantum}]
\begin{itemize}
\item[(i)] There exists a unique collection $Y_R \in \T$ indexed by all reachable rigid
objects $R$ of $\CC_w$, such that:
\begin{itemize}
\item[(a)] 
if $R,S \in \add(T)$ for some (reachable) $\CC_w$-cluster-tilting module $T$, then 
\begin{equation}\label{eq11}
Y_R Y_S = q^{[R,S]} Y_{R\oplus S} = q^{[R,S]-[S,R]}Y_S Y_R
\end{equation}
\item[(b)] if $M$ and $L$ are indecomposable rigid modules related by a mutation, with corresponding
non-split short exact sequences 
$0 \to M \to E' \to L \to 0$ and $0\to L \to E'' \to M\to 0$, then
\begin{equation}\label{eq12}
 Y_{L}Y_M = q^{[L,M]}(q^{-1}Y_{E'} + Y_{E''}).
\end{equation}
\end{itemize}
\item[(ii)] The subalgebra $\A_q(\CC_w)$ of $\T$ generated by the $Y_R$'s is a quantum
cluster algebra, in the sense of Berenstein-Zelevinsky.
The $Y_R$'s are its quantum cluster monomials (up to rescaling by some powers of $q$).
\end{itemize}
\end{theorem}
Observe that (\ref{eq11}) and (\ref{eq12}) are $q$-analogues of Theorem~\ref{Th_mult} (i) and (ii),
but in contrast with the classical case where we have defined a regular function $\varphi_X$
for every object $X$ of $\CC_w$, here we only have elements $Y_M$ of the quantum cluster algebra
for the (reachable) \emph{rigid} modules $M$.

\begin{example}\label{example7.2}
We consider again the category $\CC_w$ associated with $\g = \slchap_2$ and $w=s_2s_1s_2s_1$.
For the module $V_\i$ described in Example~\ref{Example4.4}, we calculate
the matrix 
\[
L := [\l_{ij}] = 
\left(
\begin{matrix}
0&-2&-2&-4\\
2&0&0&-2\\
2&0&0&-4\\
4&2&4&0 
\end{matrix}
\right) 
\]
encoding the $q$-commutation relations of the variables $Y_i=Y_{V_i}$. 
The two mutations discussed in Example~\ref{Example4.12} give rise
to the following quantum exchange formulas:
\[
Y_{V_1}Y_{T_1} = q^{-2}Y_{V_2}^2+Y_{V_3},
\qquad
Y_{V_2}Y_{T_2} = q^{-2}Y_{T_1}^2+Y_{V_4}. 
\]
\end{example}

\subsection{Quantum coordinate rings}\label{sect7}

Let $U_q(\g)$ be the Drinfeld-Jimbo quantized enveloping algebra of $\g$
over $\Q(q)$, with its subalgebra $U_q(\n)$. 
For a dominant weight $\l$, and $u,v\in W$ such that $u(\l)\le v(\l)$
for the Bruhat ordering, we have introduced \cite[\S5.2]{GLSquantum}
a \emph{unipotent quantum minor}
\[
D^q_{u(\l),v(\l)}\in U_q(\n).
\]
Following Lusztig, and De Concini, Kac, and Procesi, there is a 
well-defined subalgebra $U_q(\n(w))$ of $U_q(\n)$ generated 
by certain quantum root vectors $E(\b)$ labelled by the roots
$\b\in\De_w$. 
As suggested by the notation, this is a $q$-analogue of
the enveloping algebra of the nilpotent algebra $\n(w)$.
Putting 
\[
 \l_k = s_{i_1}\cdots s_{i_k}(\varpi_{i_k}),\qquad (1\le k\le r),
\]
the element $E(\b_k)$ is in fact equal to $D^q_{\l_{k^-},\l_k}$
up to a scaling factor (\cite[Proposition 7.4]{GLSquantum}).
Thus $U_q(\n(w))$ can also be regarded as a $q$-analogue of
the coordinate ring $\C[N(w)]$ of the unipotent group $N(w)$.
We can now state a quantum analogue of Theorem~\ref{thm5.3}: 
\begin{theorem}[\cite{GLSquantum}]\label{theoquantum}
Let $M_\i = M_1\oplus \cdots \oplus M_r$ be as in \S\ref{section4.7}.
The assignment 
\[
Y_{M_k} \mapsto D^q_{\l_{k^-},\l_k},\qquad (1\le k\le r),
\]
extends to an
algebra isomorphism $\iota$ from the quantum cluster algebra $\A_q(\CC_w)$
to the quantum coordinate ring $U_q(\n(w))$. 
Moreover we have
\[
\iota(Y_{V_k}) = D^q_{\varpi_{i_k},\l_k},\qquad (1\le k\le r). 
\]
\end{theorem}
The proof uses again in a crucial manner the
explicit sequence of mutations from $V_\i$ to $T_\i$ given by
Theorem~\ref{Theorem4.11}. The corresponding sequence of 
mutations for the unipotent quantum minors is similar to a 
$q$-deformation of a $T$-system, like those appearing in the
representation theory of quantum affine algebras (see \cite{KNS}).

When $\g$ is a simple Lie algebra of type $A, D, E$, and $w=w_0$
the longest element of $W$, Theorem~\ref{theoquantum} shows that $U_q(\n) = U_q(\n(w_0))$
has the structure of a quantum cluster algebra. More generally
if $w=w_0w_0^K$, where $w_0^K$ is the longest element of some
parabolic subgroup $W_K$ of $W$, then $U_q(\n(w_0w_0^K))$ 
can be regarded as the quantum coordinate ring of a big cell in
the corresponding partial flag variety, and by Theorem~\ref{theoquantum}
it also carries a quantum cluster algebra structure. 

Theorem~\ref{theoquantum} can be seen as an important step towards 
a conjecture of Berenstein and Zelevinsky \cite[Conjecture 10.10]{BFZ},
which gives a candidate for a quantum cluster structure on the quantum
coordinate ring of any double Bruhat cell $G^{v,w}$ for an arbitrary semisimple
group $G$ (not necessarily simply-laced). 
However one should pay attention to the fact that the relation
between $G^{e,w}$ and its reduced counterpart $L^{e,w}=N^w$ is not as
straightforward in the quantum case as it was in the classical case.
Indeed, to pass from the quantum coordinate ring of $N^w$
to that of $G^{e,w}$, we need to replace our \emph{unipotent} quantum
minors $D^q_{u(\l),v(\l)}$ by  ordinary ones $\De^q_{u(\l),v(\l)}$, 
which satisfy slightly different commutation relations
(see \cite[\S11.1]{GLSquantum}).

Before \cite{GLSquantum}, only a few examples 
of ``concrete'' quantum cluster algebras had appeared in the literature.
Grabowski and Launois \cite{GL} showed that the quantum coordinate rings of
the Grassmannians $\Gr(2,n) \ (n\ge 2)$, $\Gr(3,6)$, $\Gr(3,7)$, and $\Gr(3,8)$
have a quantum cluster algebra structure.
Lampe \cite{La,La2} proved two particular instances of Theorem~\ref{theoquantum}, namely
when $\g$ has type $A_n$ or $A_1^{(1)}$ and $w=c^2$ is the square of a Coxeter element. 
The existence of a quantum cluster structure on every algebra $U_q(\n(w))$ in the
general Kac-Moody case had been conjectured by Kimura \cite[Conj.1.1]{Ki}.

\begin{example}\label{exam7.4}
As in Example~\ref{example7.2}, we consider the quantum
cluster algebra $\A_q(\CC_w)$ for $\g=\slchap_2$
and $w=s_2s_1s_2s_1$. It was shown by Lampe \cite{La}
that in this case the images under $\iota$ of all the quantum cluster variables 
belong to Lusztig's dual canonical basis of $U_q(\n(w))$. 
\end{example}

Generalizing this example, Kimura and Qin \cite{KQ} have recently shown that when $w$ is 
the square of a Coxeter element and $\g$ is an arbitrary symmetric Kac-Moody algebra, 
the quantum cluster monomials of $U_q(\n(w))$ belong to the dual canonical basis.


\section{Related topics}

In this final section, we briefly mention a few additional
topics not covered in this survey, and refer the reader to the relevant references.

\subsection{Partial flag varieties}

Let us assume that $\g$ is finite-dimensional, so that $G$ is a
complex simple and simply-connected algebraic group. 
Let $B_-$ denote the Borel subgroup with unipotent radical $N_-$.
We fix a non-empty subset $J$ of $I$ and
we denote its complement by $K=I\setminus J$.
Let $B_-^K$ be the standard parabolic 
subgroup of $G$ generated by $B_-$
and the one-parameter subgroups 
\[
x_k(t),\qquad (k\in K,\,t\in\C).
\]
Consider the \emph{partial flag variety} $B_-^K\backslash G$.
For example, when $\g$ is of type $A_n$ and $J= \{j\}$, 
then $B_-^K\backslash G$ is the Grassmannian variety parametrizing $j$-dimensional
subspaces of $\C^{n+1}$.
When $\g$ is of type $D_n$ and $J= \{n\}$, 
then $B_-^K\backslash G$ is a smooth quadric in $\P^{2n-1}(\C)$.

With our assumption, $\L$ is a preprojective algebra of Dynkin type,
hence $\L$ is a finite-dimensional Frobenius algebra, and $\nil(\L) = \md(\L)$.
Let $Q_j$ be the indecomposable injective $\L$-module with socle $S_j$.
Put $Q_J := \oplus_{j\in J} Q_j$, and denote by $\Sub(Q_J)$
the full subcategory of $\md(\L)$ whose objects are submodules of
direct sums of finitely many copies of $Q_J$.
This is a Frobenius stably 2-Calabi-Yau category.

In analogy with \S\ref{sect5}, it is shown in \cite{GLSflag} that $\Sub(Q_J)$ provides a categorical
model for a cluster algebra structure on the coordinate ring of an open
cell of $B_-^K\backslash G$ (see also \cite[\S17]{GLSKM}, \cite{BIRS}). 
Moreover, by suitably extending the coefficient ring, we can lift it 
to a cluster algebra structure on the multi-homogeneous coordinate
ring of any  type $A$ partial flag variety.
This generalizes previous results of Gekhtman, Shapiro, Vainshtein \cite{GSV},
of Scott \cite{Sc0} (for type $A$ Grassmannians), and of Berenstein, Fomin, 
Zelevinsky \cite{BFZ} (for complete flag varieties $B_-\backslash G$).

\subsection{Total positivity}

As already mentioned, Lusztig's theory of total positivity for real algebraic
groups \cite{L2, L3} 
is one of the initial motivations for introducing cluster algebras
\cite{FominSurvey, FZSurvey}. 
The basic idea is that if $X$ is a variety having a totally positive part $X_{>0}$
in Lusztig's sense, then the coordinate ring of $X$ should have a cluster structure
such that each cluster gives rise to a positive coordinate system.

For instance, this has been verified for double Bruhat cells $G^{u,v}$ 
in \cite{BFZ} (see \cite[Remark 4.8]{FZSurvey}).
As a consequence, it holds for complete flag varieties $B_-\backslash G$.
It was also proved for Grassmannians in \cite{Sc0}.
Recently, using the results of \cite{GLSflag, GLSKM}, Chevalier \cite{C}
has extended this result to the partial flag varieties $B_-^K\backslash G$
in simply-laced type.

\subsection{Canonical bases}

The precise relation between cluster algebras and dual canonical bases coming from
the theory of quantum groups is still elusive, and remains a subject of active  
research. It is expected that when an algebra has both a cluster structure
and a dual canonical basis, like the coordinate ring $\C[N]$ of a maximal
unipotent subgroup of a semisimple group, then the cluster monomials should
form a subset of the dual canonical basis \cite{FZ1}.

In our setting, Theorem~\ref{Theorem5.1} shows that all cluster
monomials belong to the dual \emph{semi}canonical basis of $\C[N(w)]$. 
But it is known that in general the canonical basis differs from the 
semicanonical one \cite{GLS1}. To prove that the conjecture of Fomin and Zelevinsky 
holds in this case, one would have to understand better the intersection
of these two bases. In this direction, we have formulated the 
\emph{open orbit conjecture}. 
Let $Z$ be an irreducible component of $\L_\d$. The group 
$\GL_\d := \prod_{i\in I} GL(d_i,\C)$ acts naturally on $\L_\d$
(its orbits are in natural one-to-one correspondence with isoclasses
of $\L$-modules of dimension vector $\d$).

\begin{conjecture}[\cite{GLSKM}]
Let $Z$ be an irreducible component of $\L_\d$,
and let $\varphi_Z$ be the associated dual 
semicanonical basis vector.
If $Z$ contains an open $\GL_\d$-orbit, then
$\varphi_Z$ belongs to the dual canonical basis of $\C[N]$.
\end{conjecture}

Here, by dual canonical basis we mean the specialization at $q=1$
of the dual canonical basis of $U_q(\n)$.
In view of \S\ref{sect7}, one may also conjecture
that the quantum cluster monomials of $U_q(\n(w))$
belong to the dual canonical basis of $U_q(\n)$
(see \eg\,Example~\ref{exam7.4}).

\subsection{Grothendieck rings of quantum affine algebras}

Let us assume again that $\g$ is finite-dimensional.
Let $L\g = \g \otimes \C[t,t^{-1}]$ be the loop algebra of $\g$,
and let $U_q(L\g)$ denote the quantum analogue of its enveloping
algebra, introduced by Drinfeld and Jimbo. 
Here we assume that $q\in\C^*$ is not a root of unity.

In \cite{HL}, Hernandez and Leclerc have introduced a cluster algebra structure on the Grothendieck
ring of a certain tensor category $\CC_1$ of finite-dimensional 
$U_q(L\g)$-modules. This cluster algebra
has finitely many cluster variables, and its cluster type
(in the sense of Theorem~\ref{thFZ}~(ii)) coincides with
the Lie type of $\g$. 
Moreover, the cluster monomials coincide with the classes of simple objects
of $\CC_1$.  
This was shown in \cite{HL} for type $A_n$ and $D_4$,
and later extended to all $A, D, E$ types by Nakajima \cite{N}.
It follows that the fusion rules of $\CC_1$ are completely
encoded by the combinatorics of cluster variables.

It is conjectured in \cite{HL} that similar results hold
for more general tensor categories $\CC_{\ell}$ parametrized 
by arbitrary integers $\ell >0$. Moreover, it follows from
\cite{GLSquantum, HL2} that the quantum Grothendieck rings of 
the categories $\CC_\ell\ (\ell \le h/2-1)$, where $h$ is the Coxeter
number of $\g$, have the structure
of a quantum cluster algebra.


\frenchspacing
\begin{thebibliography}{7}

\bibitem{BFZ0}
A. Berenstein, S. Fomin and A. Zelevinsky,
\emph{Parametrizations of canonical bases and totally positive matrices},
Adv. Math. 122 (1996), 49--149.

\bibitem{BFZ}A. Berenstein, S. Fomin, A. Zelevinsky,
\emph{Cluster algebras III: Upper bounds and double Bruhat cells},
Duke Math. J. {\bf 126} (2005), 1--52.

\bibitem{BZ0}
A. Berenstein, A. Zelevinsky,
\emph{String bases for quantum groups of type $A_r$},
in I.M. Gelfand seminar, Adv. Soviet Math. {\bf 16}, pp. 51--89,
AMS 1993.

\bibitem{BZ}
A. Berenstein, A. Zelevinsky,
{\em Quantum cluster algebras},
Adv. Math.  {\bf 195}  (2005),  405--455.

\bibitem{BIRS}A. Buan, O. Iyama, I. Reiten, J. Scott,
\emph{Cluster structures for $2$-Calabi-Yau categories and unipotent
groups}, Compos. Math. {\bf 145} (2009), 1035--1079.

\bibitem{BIRS2}
A. Buan, O. Iyama, I. Reiten, D. Smith, 
\emph{Mutation of cluster-tilting objects and potentials},
Amer. J. Math. {\bf 133} (2011), 835--887. 


\bibitem{C}
N. Chevalier, 
\emph{Total positivity criteria for partial flag varieties},
J. Algebra {\bf 348} (2011), 402--415.

\bibitem{D}
L. Demonet,
\emph{Categorification of skew-symmetrizable cluster algebras},
 Algebr. Represent. Theory  {\bf 14}  (2011),  1087--1162.

\bibitem{DWZ1}
H. Derksen, J. Weyman, A. Zelevinsky, 
\emph{Quivers with potentials and their representations I: mutations},
Selecta Math. {\bf 14} (2008), 59--119. 

\bibitem{DWZ2}
H. Derksen, J. Weyman, A. Zelevinsky, 
\emph{Quivers with potentials and their representations II: applications to cluster algebras},
J. Amer. Math. Soc.  {\bf 23}  (2010), 749--790. 



\bibitem{FominSurvey}
S. Fomin,
\emph{Total positivity and cluster algebras},
Proceedings of the International Congress of Mathematicians,
Volume II, Hindustan Book Agency, 2010, 125--145.

\bibitem{FZTP}
S. Fomin, A. Zelevinsky, 
\emph{Double Bruhat cells and total positivity},
J. Amer. Math. Soc. 12 (1999), 335--380.

\bibitem{FZ1}S. Fomin, A. Zelevinsky, 
\emph{Cluster algebras I: Foundations},
J. Amer. Math. Soc. {\bf 15} (2002), 497--529.

\bibitem{FZ2}S. Fomin, A. Zelevinsky, 
\emph{Cluster algebras II: Finite type classification},
Invent. Math. {\bf 154} (2003), 63--121.

\bibitem{FZ4}S. Fomin, A. Zelevinsky,
\emph{Cluster algebras IV: Coefficients},
Compos. Math.  {\bf 143}  (2007),  112--164.

\bibitem{FZSurvey} S. Fomin, A. Zelevinsky,
\emph{Cluster algebras: notes for the CDM-$03$ conference},
Current developments in mathematics, 2003, Int. Press, Somerville,
MA 2003, 1--34. 

\bibitem{GSV} M. Gehktman, M. Shapiro, A. Vainshtein,
\emph{Cluster algebras and Poisson geometry},
Moscow Math. J. {\bf 3} (2003), 899--934.

\bibitem{GSW} M. Gehktman, M. Shapiro, A. Vainshtein,
\emph{Cluster algebras and Poisson geometry},
Mathematical surveys and monographs {\bf 167},
AMS 2010.

\bibitem{GLS1}C. Geiss, B. Leclerc and J. Schr\"oer, 
\emph{Semicanonical bases and preprojective algebras},
Ann. Scient. \'Ec. Norm. Sup. 
\textbf{38} (2005), 193--253.

\bibitem{GLS2}C. Geiss, B. Leclerc and J. Schr\"oer, 
\emph{Rigid modules over preprojective algebras},
Invent. Math., \textbf{165} (2006), 589--632.

\bibitem{GLS3}C. Geiss, B. Leclerc, J. Schr\"oer, 
\emph{Auslander algebras and initial seeds for cluster algebras},
J. London Math. Soc., {\bf 75} (2007), 718--740.

\bibitem{GLS4}C. Geiss, B. Leclerc, J. Schr\"oer, 
\emph{Semicanonical bases and preprojective algebras II: A
multiplication formula},
Compositio Math., {\bf 143} (2007), 1313--1334.

\bibitem{GLSUni1}
C. Geiss, B. Leclerc and J. Schr\"oer,
\emph{Cluster algebra structures and semicanonical bases for unipotent
groups}, 
arXiv:math/0703039 (2007), 121p.

\bibitem{GLSflag}
C. Geiss, B. Leclerc, J. Schr\"oer,
{\em Partial flag varieties and preprojective algebras},
Ann. Inst. Fourier (Grenoble) {\bf 58} (2008), 825--876.

\bibitem{GLSKM}
C. Geiss, B. Leclerc, J. Schr\"oer, 
\emph{Kac-Moody groups and cluster algebras},
Adv. Math., {\bf 228} (2011), 329--433.

\bibitem{GLSCA}
C. Geiss, B. Leclerc, J. Schr\"oer, 
\emph{Generic bases of cluster algebras and the Chamber Ansatz},
J. Amer. Math. Soc.  {\bf 25}  (2012), 21--76.


\bibitem{GLSquantum}
C. Geiss, B. Leclerc, J. Schr\"oer, 
\emph{Cluster structures on quantum coordinate rings},
arXiv:1104.0531 (2011).

\bibitem{GP}I. M. Gelfand, V.A. Ponomarev,
\emph{Model algebras and representations of graphs},
Functional Anal. Appl. {\bf 13} (1979), 157--165.

\bibitem{GL}
J. Grabowski, S. Launois,
{\em Quantum cluster algebra structures on quantum Grassmannians and their
quantum Schubert cells: the finite-type cases},
Int. Math. Res. Not. IMRN  {\bf 10} (2011), 2230--2262.

\bibitem{H2}
D. Happel,
\emph{Triangulated categories in the representation theory of 
finite-dimensional algebras},
London Mathematical Society Lecture Note Series, 119. 
Cambridge University Press, Cambridge, 1988. x+208 pp.


\bibitem{HL} {D. Hernandez, B. Leclerc}, 
{\em Cluster algebras and quantum affine algebras},
Duke Math. J. \textbf{154} (2010), 265--341.

\bibitem{HL2} {D. Hernandez, B. Leclerc}, 
\emph{Quantum Grothendieck rings and derived Hall algebras},
arXiv:1109.0862 (2011).

\bibitem{KP}
V. Kac and D. Peterson,
\emph{Regular functions on certain infinite-dimensional groups},
Arithmetic and geometry, Vol II, 141--166, Progr. Math., 36,
Birkh\"auser Boston, Boston, MA, 1983.

\bibitem{Keller} B. Keller,
\emph{Cluster algebras and derived categories},
arXiv:1202.4161, (2012), 60 p.

\bibitem{Ki}
Y. Kimura,
{\em Quantum unipotent subgroup and dual canonical basis},
Kyoto J. Math. {\bf 52} (2012), 277--331. 

\bibitem{KQ}
Y. Kimura, F. Qin,
{\em Graded quiver varieties, quantum cluster algebras and dual canonical basis},
arXiv:1205.2066.

\bibitem{Ku}
S. Kumar,
\emph{Kac-Moody groups, their flag varieties and representation theory},
Progress in Mathematics, 204. Birkh\"auser Boston, Inc., Boston, MA, 2002. 

\bibitem{KNS}
A. Kuniba, T. Nakanishi, J. Suzuki,
\emph{$T$-systems and $Y$-systems in integrable systems},
J. Phys. A  {\bf 44}  (2011), 146 p. 

\bibitem{La}
P. Lampe,
{\em A quantum cluster algebra of Kronecker type and the dual canonical
basis},
Int. Math. Res. Not. IMRN  {\bf 13} (2011), 2970--3005. 

\bibitem{La2}
P. Lampe,
{\em Quantum cluster algebras of type A and the dual canonical basis},
arXiv:1101.0580, (2010), 46 p. 

\bibitem{Lu0}G. Lusztig,
\emph{Quivers, perverse sheaves, and quantized enveloping algebras}, 
J. Amer. Math. Soc. {\bf 4} (1991), 365--421.



\bibitem{L2} G. Lusztig, 
\emph{Total positivity in reductive groups},  
In \emph{Lie theory and geometry}, 
Progr. Math., \textbf{123}, Birkh\"auser Boston, 1994, 531--568.

\bibitem{L3} G. Lusztig, 
\emph{Total positivity in partial flag manifolds}, 
Represent. Theory \textbf{2} (1998), 70--78.


\bibitem{Lu1}G. Lusztig,
\emph{Semicanonical bases arising from enveloping algebras}, 
Adv. Math. {\bf 151} (2000),  129--139.


\bibitem{N}H. Nakajima,
\emph{Quiver varieties and cluster algebras},
Kyoto J. Math.  51  (2011),  71--126. 

\bibitem{Ri}C.M. Ringel,
\emph{The preprojective algebra of a quiver},  Algebras
and modules, II (Geiranger, 1996),  467--480, CMS Conf. Proc., {\bf 24},
Amer. Math. Soc., Providence, RI, 1998.

\bibitem{Sc0} J. Scott,
\emph{Grassmannians and cluster algebras},
Proc. London Math. Soc. {\bf 92} (2006), 345--380.

\bibitem{Sc} J. Scott,
\emph{Block-Toeplitz determinants, chess tableaux, and the type $\hat{A_1}$ 
Geiss-Leclerc-Schr\"oer $\varphi$-map},
 arXiv:0707.3046, (2007), 23 p.

\end{thebibliography}
\end{document}